\documentclass{amsart}
\usepackage{amssymb,amsmath,amscd,xy}
\newtheorem{theorem}{Theorem}[section]
\newtheorem{lemma}[theorem]{Lemma}
\newtheorem{corollary}[theorem]{Corollary}
\newtheorem{definition}[theorem]{Definition}

\newtheorem{remark}[theorem]{\it Remark}
\newtheorem{example}[theorem]{Example}
\newtheorem{proposition}[theorem]{Proposition}
\newtheorem{maintheorem}[theorem]{Main Theorem}

\xyoption{arrow}

\xyoption{matrix}

\newcommand{\tableau}[2]{\begin{array}{|c|}\hline #1 \\ \hline #2 \\ \hline \end{array}}

\def\gr{\mathrm{gr}}
\def\Gr{\mathrm{Gr}}
\def\GL{\mathrm{GL}}
\def\PGL{\mathrm{PGL}}
\def\SL{\mathrm{SL}}
\def\CP{\mathbb{C}P}
\def\PSL{\mathrm{PSL}}
\def\SO{\mathrm{SO}}
\def\wt{\mathbf{wt}}
\def\t{\mathrm}

\def\C{\mathbb{C}}
\def\E{\mathbb{E}}
\def\R{\mathbb{R}}
\def\Z{\mathbb{Z}}
\def\N{\mathbb{N}}
\def\T{\mathbb{T}}

\def\CP{\mathbb{CP}}

\def\br{\mathbf{r}}
\def\bs{\mathbf{s}}
\def\bd{\mathbf{d}}
\def\be{\mathbf{e}}
\def\ba{\mathbf{a}}
\def\bb{\mathbf{b}}
\def\bk{\mathbf{k}}
\def\bt{\mathbf{t}}
\def\bc{\mathbf{c}}
\def\bx{\mathbf{x}}

\def\q{/\!/}

\def\cq{/\!/}
\def\Proj{\mathrm{Proj} \,}

\def\Mr{M_\mathbf{r}(\R^2)}

\begin{document}

\title[The projective invariants of ordered points on the line]
{The projective invariants of ordered points on the line }

\author{Benjamin Howard, John Millson, Andrew Snowden and Ravi Vakil*}
\thanks{*The first two authors were partially supported by 
NSF grant DMS-0104006 and the last author
was partially supported by NSF Grant DMS-0228011}
\date{February 11, 2006.}

\maketitle

\begin{abstract}
  The space of $n$ (ordered) points on the projective line, modulo
  automorphisms of the line, is one of the most important and
  classical examples of an invariant theory quotient, and is one of
  the first examples given in any course.  Generators for the ring of
  invariants have been known since the end of the nineteenth century,
  but the question of the relations has remained surprisingly open,
  and it was not even known that the relations have bounded degree.
  We show that the ideal of relations is generated in degree at most
  $4$, and give an explicit description of the generators.  The result
  holds for arbitrary weighting of the points.  If all the weights are
  even (e.g.\ in the case of equal weight for odd $n$), we show that
  the ideal of relations is generated by quadrics.  The proof is by
  degenerating the moduli space to a toric variety, and following an
  enlarged set of generators through this degeneration.  

In a later
  work \cite{HowardMillsonSnowdenVakil}, by different means, we will
  show that the projective scheme is cut out by certain natural
  quadrics unless $n=6$ and each point has weight one.
\end{abstract}

\tableofcontents

\section{Introduction}\label{introduction}

The space of $n$ (ordered) points on the projective line, modulo
automorphisms of the line, is one of the most important and classical
examples of an invariant theory quotient, and is one of the first
examples given in any course (see \cite[\S 2]{ms}, \cite[\S 3]{git},
\cite[Ch.~11]{Dolgachev}, \cite[Ch.~I]{DO}, \dots). 

The generators have been known for a long time: in 1894 Kempe, see
\cite{Kempe}, proved that for the case of $n$ points,
the lowest degree invariants generate the ring of
invariants. This paper was said by Howe in \cite{Howe}, pg.\ 156,
to be the ``deepest result'' of classical invariant theory (in the
sense that it is the only result that Howe could not prove from
standard constructions in representation theory).
However, the question of the relations has remained surprisingly open.
It was not even known that the relations have bounded degree as $n$
varies.

In this context, it is most natural to consider more generally the
space of {\it weighted} points on $\CP^1$.  Let the $i$th point be
weighted by $r_i$ and let $\br = (r_1,\ldots,r_n) \in (\Z^+)^n$. Let
$R_\br$ denote the ring of projective invariants.  The weights can be
interpreted as parametrizing the very ample line bundles of
$(\CP^1)^n$. The space of weighted points (for real weights) was
considered by Deligne and Mostow \cite{DeligneMostow}. In their theory the weights are
the parameters for the hypergeometric function. We will see below
another reason (coming from Euclidean geometry) why introducing
weights is natural --- the weights are the side-lengths of polygonal
linkages in Euclidean three-space.

We first observe that Kempe's Theorem (that the ring of invariants is generated
in ``lowest degree'') remains true for all integral
weights (Theorem~\ref{generatorstheorem}), incidentally giving a short, possibly new, proof of Kempe's Theorem.
We dub these generators the {\em Kempe generators} of the ring.
Our main result on the {\it
relations} of the invariants is the following.

\begin{maintheorem}  For any weights $\br = (r_1,\ldots,r_n)$, the ideal
of relations in the ring of invariants is generated by relations of degree at
most four.  If all the weights are even, then the ideal of relations is generated
in degree two.
\end{maintheorem}

The proof is completed in \S \ref{ring} (see Theorems~\ref{relations}
and \ref{maintheorem}).  We give an explicit description of a
generating set (suitable for example for computation) in \S
\ref{lifting}, and examples at the end of \S \ref{ring}.

The heart of the paper is \S \ref{toricring} and \S \ref{ring}.  In
these sections we enlarge the set of generators in order that the
leading terms of this larger set (relative to the filtration by
Lakshmibai-Gonciulea degree defined in \S \ref{toric_degeneration}) is
a set of generators for the associated graded ring (the homogeneous
coordinate ring of the toric fiber of our toric degeneration). The
extra generators are quadratic expressions in the Kempe
generators. Our main result proves that the ideal of relations among
this larger set of generators is generated by quadratic
relations. When these quadratic relations are rewritten in terms of
the Kempe generators we obtain relations of degrees two, three and
four. We do this by finding a normal form for monomials in the Kempe
generators.

The polytope corresponding to the toric variety has three
interpretations, which we discuss in \S \ref{threepolytopes}.  One
standard interpretation is in terms of Gel'fand-Tsetlin patterns, and another
is in terms of semistandard tableaux, but our proof uses a 
third interpretation in terms of possible lengths of diagonals emanating 
from an initial 
vertex of an $n$-gon with side lengths $\br = (r_1, \ldots, r_n)$.  As a
consequence of this polytope interpretation we get effective
bounds on the degrees of these moduli spaces, enabling us to show that
the ring of projective invariants for $n$ equally weighted points is not a complete intersection unless
$n \leq 4$ or $n = 6$ (\S \ref{application}).  This result seems not
to be geometrically obvious.

In our next paper \cite{HowardMillsonSnowdenVakil}, we prove that the
moduli space is in fact cut out by combinatorially/geometrically
obvious quadratic relations, with the single exception of $\br
= (1,1,1,1,1,1)$ (corresponding to the Segre cubic threefold).
However, we could not prove that these quadratic relations generate
the {\em ideal} of relations.

Thus the understanding of the ring of projective invariants of ordered
points on the projective line is now quite satisfactory.  This is in
contrast with the equally classical (and much more complicated)
question of {\it unordered} points.  (This is equivalent to
considering homogeneous degree $n$ polynomials in two variables ---
the subject of ``binary quantics,'' see \cite{KungRota}.)
 Mumford has written: ``This is an
amazingly difficult job, and complete success was achieved only for $n
\leq 6$'' \cite[pg.~77]{git}.  Later, ``by an extraordinary tour de
force'' (Mumford, loc.\ cit.) Shioda \cite{Shioda} dealt with the case
$n=8$.  In the unordered case, a generating set for general $n$ is not
known (and not even known to have bounded degree as $n$ varies).

The case of unordered points with even $n$ essentially corresponds to
the ring of hyperelliptic modular forms of genus $(n-2)/2$, and their
relations.  Our case of ordered points, where $n$ is even and the the
weights are even, essentially corresponds to the ring of hyperelliptic
level two modular forms, and this paper completely describes
generators of the ideal of relations among these forms.

\vskip 12pt

\noindent {\em Acknowledgments.}
Yi Hu has  pointed out that our ``side-splitting map'' (see \S 3.5)
is a special case of the splitting map introduced in \S 2.3 of
\cite{Hu}. In particular our Proposition 3.5 is a special case
of his Proposition 2.11.
We would like to thank Philip Foth for many helpful comments about
toric degenerations. This project was inspired by reading
\cite{FothHu}. In \cite{HowardManonMillson} the authors will prove the
conjecture of \cite{FothHu} which describes the toric fiber of the
toric degeneration of the moduli space of points on the line (as a
topological space) as the collapsed space introduced in
\cite{KamiyamaYoshida}.  The fourth author thanks Allen Knutson and
Diane Maclagan for enlightening conversations. Also we thank
Allen Knutson for pointing out that the first toric degenerations
were given in Proposition 11.10 (page 104) of \cite{Sturmfels}.
Finally we would like to thank Roger Howe for telling us about \cite{Howe}.

\section{G.I.T. quotients}\label{mumford_quotient}

\subsection{Definition of G.I.T. quotient}\label{Mumford_Defn}
We refer the reader to \cite{Dolgachev} for additional details.
Suppose that $G$ is a reductive algebraic group, $V$ is a
quasi-projective variety, and $\eta : G \times V \to V$ is regular
action of $G$. Let $\pi : \mathcal{L} \to V$ be an ample line
bundle over $V$. A $G$--linearization of $\mathcal{L}$ is a
regular action $\widetilde{\eta} : G \times \mathcal{L} \to
\mathcal{L}$ which is linear on fibers and makes the following
diagram commute:
\[
\begin{CD}
   G \times \mathcal{L} @>\widetilde{\eta}>> \mathcal{L} \\
   @V{id \times \pi}VV                 @VV{\pi}V \\
   G \times V @>\eta>> V
\end{CD}
\]

Given such a linearization, we automatically get linearizations on
all tensor powers $\mathcal{L}^{\otimes N}$ of $\mathcal{L}$. Thus
$G$ has an action on sections $s$ of $\mathcal{L}^{\otimes N}$
given by $$(g \cdot s)(x) = g \cdot s(g^{-1} \cdot x) =
\widetilde{\eta}(g,s(\eta(g^{-1},x))).$$ Let
$\Gamma(V,\mathcal{L}^{\otimes N})^G$ denote the $G$--invariant
sections of $\mathcal{L}^{\otimes N}$. The G.I.T. quotient $V
\q_{\widetilde{\eta}}\: G$ is defined as
$$V \q_{\widetilde{\eta}}\: G = \t{Proj} \Big(
\bigoplus_{N=0}^\infty \Gamma(V,\mathcal{L}^{\otimes N})^G
\Big).$$  If the linearization is understood, sometimes we denote
$V \q_{\widetilde{\eta}} G$ by $V \q G$.

If $s$ is a section let $\t{supp}(s) = \{x \in V \mid s(x) \neq
0\}$. The set $V_{\widetilde{\eta}}^{ss}$ of semistable points of
$V$ is defined as
$$V_{\widetilde{\eta}}^{ss} = \bigcup_{N \geq 0}
\bigcup_{s \in \Gamma(V,\mathcal{L}^{\otimes
N})^G}
 \t{supp}(s).$$ 
(To be more precise, we need the distinguished open 
subset $\t{supp}(s)$ corresponding to $s$ 
to be affine, but this will be true in the two cases
relevant for us, when $X$ is affine or projective.)
If $x$ is a semistable point let $\overline{G \cdot x}$ be the
(Zariski) closure of the orbit $G \cdot x$ in
$V_{\widetilde{\eta}}^{ss}$. As a topological space $V
\q_{\widetilde{\eta}}\: G$ is the quotient space of
$V_{\widetilde{\eta}}^{ss}$ where points $x,y$ are identified iff
$\overline{G \cdot x}$ and $\overline{G \cdot y}$ intersect nontrivially.

\subsection{The Gel'fand-MacPherson
correspondence}\label{GMcorrespondence}

The Gel'fand--MacPherson correspondence says that a G.I.T. quotient of
Grassmannian space $\Gr_k(\C^n)$ by the torus $(\C^\ast)^n$ is
isomorphic to a G.I.T. quotient of the product space $(\CP^{k-1})^n$ by
the diagonal action of $\t{PGL}(k,\C)$.

Let $\mathcal{L}$ be the trivial line bundle $\C^{n \times k}
\times \C \rightarrow \C^{n \times k}$.  Given any group $G$
acting on $\C^{n \times k}$, a character $\chi:G \to \C^\ast$
defines a linearization of $\mathcal{L}$ by $g \cdot (A,z) = (g
\cdot A, \chi(g) z)$.

The group $\GL(k,\C)$ acts on the right of $\C^{n \times k}$ by
matrix multiplication. The group $T$ of diagonal matrices in
$\GL(n)$ acts on the left of $\C^{n \times k}$. Let $\det^a :
\GL(k,\C) \rightarrow \C^\ast$ be $\det^a(g) = (\det(g))^a$ and let
$\chi_\br : T \rightarrow \C^\ast$ be
$$\chi_\br(\t{diag}(z_1,\ldots,z_n)) = \prod_{i=1}^n z_i^{r_i}$$
where $\br = (r_1,\ldots,r_n) \in \Z_+^n$. The one--dimensional
subgroup $K = \{(zI_n,z^{-1}I_k) : z \in \C^\ast\}$ of $T \times
\GL(k,\C)$ acts trivially on $\C^{n \times k}$. Let $G$ be the
quotient of $T \times \GL(k,\C)$ by $K$. The character $\chi_\br
\times \det^a$ descends to $G$ iff $|\br|= \sum_i r_i = ka$, and we
assume that is the case so that we have a $G$--linearization of
the trivial line bundle.

Let $\mathcal{L}_{k,n}$ be the ample generator of the Picard group
of $\Gr_k(\C^n)$; we may realize the total space of
$\mathcal{L}_{k,n}$ by equivalence classes $V^{n \times k} \times
\C / \sim$ where $V^{n \times k}$ is the open subset of $\C^{n
\times k}$ of matrices with independent columns and $(A,z) \sim
(Ag,\det(g)z)$ for $g \in \GL(k,\C)$. Denote the equivalence class
of $(A,z)$ by $[A,z]$.  The character $\chi_\br$ defines a
$T$--linearization of $\mathcal{L}_{k,n}^a =
\mathcal{L}_{k,n}^{\otimes a}$ by
$$t \cdot [A,z] = [tA, \chi_\br(t)z].$$

Let $\mathcal{H}$ be the ample generator of the Picard group of
$\CP^{k-1}$, and let $\mathcal{H}^\br$ be the ample line bundle over
the product $(\CP ^{k-1})^n$ given by
$$\mathcal{H}^\br = \mathcal{H}^{\otimes r_1} \boxtimes \cdots
\boxtimes \mathcal{H}^{\otimes r_n}.$$ We may identify the total
space of $\mathcal{H}$ with $(\C^k \setminus \{0\}) \times \C /
\sim$ where $(v,z) \sim (v\lambda, \lambda z)$ for $\lambda \in
\GL(1)$; let $[v,z]$ denote the equivalence class. There is a
unique linearization of $\mathcal{H}^\br$ for the (right) diagonal
action of $\PGL(k,\C)$ on $(\CP^{k-1})^n$.

By the First Fundamental Theorem of Invariant Theory
\cite[Theorem~2.1]{Dolgachev}, the homogeneous coordinate
  ring of the Grassmannian is generated by Pl\"ucker coordinates,
  hence, for any $N \geq 0$ we have
$$\Gamma(\Gr_k(\C^n),(\mathcal{L}_{k,n}^a)^{\otimes N}) \cong
\Gamma(\C^{n \times k},\mathcal{L}^{\otimes akN})^{\GL(k,\C)}$$
and consequently
$$\Gamma(\Gr_k(\C^n),(\mathcal{L}_{k,n}^a)^{\otimes N})^T \cong
\Gamma(\C^{n \times k},\mathcal{L}^{\otimes akN})^{G}.$$
On the other hand the sections of the outer tensor product of
hyperplane section bundles over $n$  copies of $\CP^{k-1}$ are products
of {\em homogeneous} polynomials in the homogeneous coordinates, that is,
$$\Gamma((\CP^{k-1})^n,(\mathcal{H}^\br)^{\otimes N}) \cong
\Gamma(\C^{n \times k},\mathcal{L}^{\otimes |\br| N})^T$$
and consequently
$$\Gamma((\CP^{k-1})^n,(\mathcal{H}^\br)^{\otimes N})^{\PGL(k,\C)} \cong
\Gamma(\C^{n \times k},\mathcal{L}^{\otimes |\br| N})^{G}.$$
Hence we have an isomorphism of the G.I.T. quotients
$$\Gr_k(\C^n) \q_{\chi_\br} T \cong \C^{n \times k} \q_{\chi_\br
\times \det^a} G \cong (\mathbb{CP}^{k-1})^n \q _{\br} \PGL(k,\C).$$

In the Introduction we took the quotient of
$(\mathbb{CP}^1)^n$ by $\SL(2,\C)$ instead of $\PGL(2,\C)$. We claim
that there exists a $\PGL(2,\C)$--linearization of $\mathcal{H}^{\br}$ if 
and only if the sum of the $r_i$'s is even. Indeed, we have a 
$\PGL(2,\C)$--linearization if and only if we have a $\PSL(2,\C)$--linearization 
if and only if we have an $\SL(2,\C)$ linearization in which $-1$ acts trivially. 
But $-1$ acts on $\mathcal{H}^{\br}$ by multiplication on each fiber by 
$(-1)^{\sum r_i}$ and the claim follows. On the other hand we have an
$\SL(2,\C)$--linearization and an $\SL(2,\C)$ quotient for any $\br$.
The two quotients coincide when they are both defined (i.e.\ the
sum of the $r_i$'s is even).

\subsection{Points on the line and Euclidean polygons}\label{moduli_space}
Here we briefly explain why the G.I.T. quotient $(\CP^1)^n \q_\br \SL(2,\C)$ 
is homeomorphic (and 
symplectomorphic in the smooth case) to the space $M_{\br}$ of congruence classes 
of Euclidean $n$-gon linkages with side-lengths $\br$. 
We refer the reader to \cite{KapovichMillson} for details. 

An $n$-gon $\be$ modulo translations in Euclidean three
space $\E^3$ is given by an $n$--tuple of vectors $\be = (e_1,
e_2, \ldots,e_n)$ (the edges) such that the following {\em closing
condition} holds
$$e_1 + e_2 + \cdots + e_n = 0.$$

It is proved in \cite{KapovichMillson} and \cite{Klyachko} that
the map $\mu(\be) = e_1 + e_2 + \cdots + e_n$ is the momentum map
for the diagonal action of $\SO(3)$ on the product $S^2(r_1) \times
\cdots \times S^2(r_n)$. Consequently the set of closed $n$-gon
linkages is the zero level set of the momentum map and after
dividing by $\SO(3)$ we obtain the symplectic quotient $S^2(r_1)
\times \cdots \times S^2(r_n)\q _{\br} \SO(3)$. It follows from a
general theorem of Kirwan-Kempf-Ness that the symplectic quotient
$S^2(r_1) \times \cdots \times S^2(r_n)\q _{\br} \SO(3)$ is
canonically homeomorphic to the G.I.T. quotient
$(\mathbb{CP}^1)^n\q_{\br} \SL(2,\C)$. A direct proof using the
conformal center of mass is given in \cite{KapovichMillson}. We
obtain
\begin{theorem}
There is a canonical homeomorphism of the space $M_{\br}$ of
$n$-gon linkages in $\E^3$ with side-lengths $\br$ modulo
orientation-preserving isometries and the moduli space
$(\mathbb{CP}^1)^n\q_{\br} \SL(2,\C)$ of weighted ordered points
on the line.
\end{theorem}


\section{The algebraic geometry of the moduli space of
polygons}\label{AlgGeomSpace}

\subsection{The Pl\"ucker embedding of the Grassmannian}
In the previous section we gave a construction of $\Gr_k(\C^n)$ as
a G.I.T. quotient of $\C^{n \times k}$.  We now look at this more
closely for $k=2$.

Recall that $\mathcal{L}$ is the trivial line bundle over $\C^{n
\times 2}$ and we have a $\GL(2,\C)$--linearization of the bundle
$\mathcal{L}^{\otimes N}$ via the character $\mathrm{det}^a$.
Specifically, for $m \in \C^{n \times 2}$, $z \in \C$ and $g \in
\GL(2,\C)$, we have
$$g \cdot (m, z)=(mg^{-1}, \det(g)^{-aN} z).$$
Given a map $\tau:\C^{n \times 2} \rightarrow \C$ we get a section of
$\mathcal{L}^{\otimes N}$ via $s_\tau(m)=(m,\tau(m))$, and all sections
are of this form.  For $g \in \GL(2,\C)$ we have
$$(g \cdot s_\tau)(m)=g \cdot s_\tau(mg)=
g \cdot (mg, \tau(mg))=(m, \det(g)^{-aN} \tau(mg)).$$ Thus the section
$s_\tau$ is invariant if and only if $\tau(mg)=\det(g)^{aN} \tau(m).$

Let $i$ and $j$ be integers, $1 \le i, j \le n$.  Define a map
$\C^{n \times 2} \rightarrow \C$ by assigning to a matrix $m$ the
determinant of the $2 \times 2$ minor of $m$ formed by taking rows
$i$ and $j$.  We denote this map with the $2 \times 1 $ tableau
$$\begin{array}{|c|} \hline i \\ \hline j \\ \hline \end{array} \, .$$
(Note that we have the identity
$$\begin{array}{|c|} \hline i \\ \hline j \\ \hline \end{array}=
- \, \begin{array}{|c|} \hline j \\ \hline i \\ \hline \end{array}
\, ,$$ so it suffices to consider tableaux with $i<j$.) The
product of $k$ of these $2 \times 1$ tableaux is denoted by a $2
\times k$ tableau:
$$\begin{array}{|c|} \hline i_1 \\ \hline j_1 \\ \hline \end{array} \cdot
\begin{array}{|c|} \hline i_2 \\ \hline j_2 \\ \hline \end{array} \cdots
\begin{array}{|c|} \hline i_k \\ \hline j_k \\ \hline \end{array} =
\begin{array}{|c|c|c|c|} \hline
i_1 & i_2 & \cdots & i_k \\ \hline j_1 & j_2 & \cdots & j_k \\
\hline
\end{array} \, .$$

Clearly if $\tau:\C^{n \times 2} \rightarrow \C$ is a map given by a
$2 \times k$ tableau then $\tau(mg)=\det(g)^k \tau(m)$; in fact any
map $\tau:\C^{n \times 2} \rightarrow \C$ satisfying $\tau(mg)=\det(g)^k
\tau(m)$ is given by a linear combination of $2 \times k$ tableau.
Equivalently, all invariant sections of
$\mathcal{L}^{\otimes N}$ (or sections of $\mathcal{H}^{\otimes
a}$, see previous section) are linear combinations of $s_\tau$ where
$\tau$ is a $2 \times aN$ tableau.

The G.I.T. quotient $\C^{n \times 2} \cq_{\det^a} \GL(2,\C)$ (which
we know to be $\Gr_2(\C^n)$) is by definition $\Proj{S_a}$ where
$S_a$ is the graded ring $\bigoplus_{N \ge 0} \Gamma(\C^{n \times
2}, \mathcal{L}^{\otimes N})^{\GL(2,\C)}$. We have a good
understanding of the generators of $S_a$:  the ring is generated
by lowest degree elements, and each graded piece is spanned by
$s_\tau$, over tableaux $\tau$ of the appropriate size.  We now describe
the relations satisfied by the $s_\tau$.  

A tableau is called \emph{semistandard} if its columns are
strictly increasing and its rows are weakly increasing.
The following two propositions are standard facts in
invariant theory, see for example \cite{Dolgachev}.  We will give additional proofs
\cite{HowardMillsonSnowdenVakil} in terms of ``Kempe graphs''.

\begin{proposition}
The relations between tableaux are generated by the
\emph{Pl\"ucker relations}:
$$\begin{array}{|c|c|} \hline a & b \\ \hline c & d \\ \hline \end{array}=
\begin{array}{|c|c|} \hline a & c \\ \hline b & d \\ \hline \end{array}+
\begin{array}{|c|c|} \hline a & b \\ \hline d & c \\ \hline \end{array}$$
where $a,b,c,d$ are any integers between 1 and $n$.  If $\tau$, $\sigma$,
and $\gamma$ are the above tableaux, then all relations in $S_1$ are
generated by the corresponding relation $s_\tau=s_\sigma+s_\gamma$.
\end{proposition}

\begin{proposition}
The $2 \times aN$ semistandard tableaux are linearly independent;
$$\{ s_\tau \mid \tau \textrm{ is $2 \times aN$, semistandard} \}$$
is a basis for the $N$th graded piece of $S_a$.
\end{proposition}

\subsection{The homogeneous ring of $(\CP^1)^n \q_\br \SL(2,\C)$}
Assume that $|\br| = \sum_i r_i = 2a$ is even.
If we quotient $\C^{n \times 2}$ by $\GL(2,\C)$ 
using the linearization $|\br|$, 
we obtain the Grassmannian
$\Gr_2(\C^n)$ with the line bundle $\mathcal{L}_{2,n}^{\otimes a}$.
Let $\br$ be a weight of the torus $T \subset \GL(n)$ and
$\chi_{\br}$ the corresponding character, which defines a $T$--linearization
of $\mathcal{L}_{2,n}^{\otimes a}$.
If we now quotient $\Gr_2(\C^n)$ by $T$ we obtain the moduli space
$(\mathbb{CP}^1)^n\q_\br \SL(2,\C)$ (see the previous section on the
Gel'fand-MacPherson correspondence.)

The quotient $\Gr_2(\C^n) \q_\br T$ is, by definition,
$\Proj{S_a^T}$, where
$S_a^T$ is the fixed subring of $S_a$ under the action of the
torus.  We now determine the invariant sections.   Let $s_\tau$ be a
section of $\mathcal{H}^{\otimes aN}$ corresponding to the tableau
$\tau$.  An element $\bt=(t_1, \ldots, t_n) \in T$ acts on $\tau$ via
$$\bt \cdot \tau = \bt^{\bc-N \br} \tau$$
where $c_i$ is the number of times $i$ occurs in $\tau$ and we use
the notation
$$\bt^{\bc}=\prod_{i=1}^n t_i^{c_i}$$
(this is most easily seen by returning to the bundle $\mathcal{L}$
over $\C^{n \times 2}$). Thus $s_\tau$ is invariant if and only if
$\bc=N \br$, {i.e.}, $c_i=N r_i$ for each $1 \le i \le n$.
For a general section $s$, write $s=\sum \alpha_i s_{\tau_i}$ where
the $\tau_i$ are semistandard. Clearly, $s$ is invariant if and only
if each $s_{\tau_i}$ is as well. Thus we have shown:

\begin{proposition}
The $N$th graded piece of the ring $S_a^T$ has for a basis $\{ s_\tau
\}$ where $\tau$ runs over the semistandard tableaux for which $i$
occurs precisely $N r_i$ times.
\end{proposition}

From now on, we are going to be working more with the tableaux
than the sections.  For convenience, we put $R_{\br}$ to be the
``ring of tableaux'' which satisfy the condition of the above
theorem.  This is clearly no loss in generality, as the map
$s:R_{\br} \rightarrow S_a^T$ which assigns to a tableau $\tau$ the
section $s_\tau$ is an isomorphism.

We also will write $\C[\Gr_2(\C^n)]$ for the homogeneous
coordinate ring of the Grassmannian with respect to the
Pl\"ucker embedding. We shall regard the elements of
$\C[\Gr_2(\C^n)]$ as tableaux with two rows and any number of
columns, with entries ranging from 1 to $n$.  The ring $R_{\br}$
is a subring of $\C[\Gr_2(\C^n)]$.

\subsection{Example}  If $\br=(1,1,\ldots,1)$ (which corresponds to
equilateral $n$-gons) a tableau is in $R_{\br}$ if and only if each
integer between 1 and $n$ occurs the same number of times in the
tableau. For example, for $n=6$ the tableau
$$\begin{array}{|c|c|c|} \hline
1 & 2 & 3 \\ \hline 4 & 5 & 6 \\ \hline
\end{array}$$
is allowed.  Note that the condition of being in $R_{\br}$ puts a
restriction on the number $k$ of columns in a tableau:  for this
choice of $\br$, if $n$ is even then $k$ must be a multiple of
$n/2$ while if $n$ odd then $k$ must be a multiple of $n$.

\subsection{The side-splitting map on rings}
We now construct a surjection of graded rings $R_{\br} \rightarrow
R_{\br'}$ for certain choices of $\br$ and $\br'$.  This map will
be called the \emph{side-splitting map} for reasons which will
become apparent later. In this section, $[x,y]$ will denote the
integers between $x$ and $y$ (inclusively).

Let $\br$ be a side-length vector of length $n$ and let $\br'$ be
a side length vector of length $n'$, where $n>n'$.  Let
$\phi:[1,n] \rightarrow [1,n']$ be a nondecreasing map such that
\begin{equation}
\label{sidesplitcond} \br_i'=\sum_{\phi(j)=i} \br_j.
\end{equation}
We define a map $\phi_*:R_{\br} \rightarrow R_{\br'}$ as follows:
let $\phi_\ast$ be the unique graded ring homomorphism such that
for any tableau $\tau \in R_{\br}$,  $\phi_*(\tau)$ is the
tableau obtained by replacing all $i$'s in $\tau$ with $\phi(i)$.  The
following proposition shows that this is well-defined.

\begin{proposition}
\label{sidesplitmap} The map $\phi_*$ is a well-defined surjection
of graded rings. Furthermore, it maps semistandard tableaux to
semistandard tableaux (or zero).
\end{proposition}

\begin{proof}
Define a map $\hat{\phi}_*:\C[\Gr_2(\C^n)] \rightarrow
\C[\Gr_2(\C^{n'})]$ by letting $\hat{\phi}_*$ act on tableaux in
the same manner that $\phi_*$ was specified to act on tableaux.
The map $\hat{\phi}_*$ will lift the map $\phi_*$ to the larger
ring $\C[\Gr_2(\C^n)]$.  The reason for doing this is that these
larger rings are simpler than the rings $R_{\br}$

Recall that the rings $\C[\Gr_2(\C^n)]$ and $\C[\Gr_2(\C^{n'})]$
are generated by $2 \times 1$ tableaux with entries in $[1,n]$
($[1,n']$ respectively) subject to the relations
$$\begin{array}{|c|} \hline a \\ \hline b \\ \hline \end{array}=
- \, \begin{array}{|c|} \hline b \\ \hline a \\ \hline \end{array}
\qquad \textrm{and} \qquad
\begin{array}{|c|c|} \hline a & b \\ \hline c & d \\ \hline \end{array}=
\begin{array}{|c|c|} \hline a & c \\ \hline b & d \\ \hline \end{array}+
\begin{array}{|c|c|} \hline a & b \\ \hline d & c \\ \hline \end{array}$$
These relations are obviously preserved by $\hat{\phi}_*$; thus
$\hat{\phi}_*$ is well-defined.

Next observe that $\hat{\phi}_*$ maps $R_{\br}$ into $R_{\br'}$.
This is an immediate consequence of the condition
\eqref{sidesplitcond} imposed on the map $\phi$.  By construction,
the restriction of the map $\hat{\phi}_*$ to the subring $R_{\br}$
is equal to the map $\phi_*$, that is to say, we have the
following commutative diagram:
$$\xymatrix{
R_{\br} \ar@{^{(}->}[r] \ar[d]_{\phi_*} & \C[\Gr_2(\C^n)]
\ar[d]^{\hat{\phi}_*} \\
R_{\br'} \ar@{^{(}->}[r] & \C[\Gr_2(\C^{n'})]. }$$ 
This proves that
the map $\phi_*$ is well-defined, since it is simply the
restriction of a well-defined map.

The fact that $\phi_*$ is surjective is obvious:  given any
tableau $\tau$ in $\C[\Gr_2(\C^{n'})]$ simply change all occurrences
of $i$ to something in $\phi^{-1}(i)$ in such a way that the
weights add up correctly.  More precisely, if $j_1, \ldots, j_k$
are the solutions to $\phi(j)=i$ then change $\br_{j_1}$ of the
$i$'s in $\tau$ to be $j_1$'s, $\br_{j_2}$ of the $i$'s to $j_2$'s, 
{etc}.  By the condition \eqref{sidesplitcond} this is guaranteed
to work. Thus $\phi_*$ is surjective.  Note, however, that
$\hat{\phi}_*$ is not necessarily surjective.

Finally, the fact that $\phi_*$ takes semistandard tableau in
$R_{\br}$ to semistandard tableau in $R_{\br'}$ (or zero) follows
immediately from the condition that $\phi$ is nondecreasing.  This
completes the proof.
\end{proof}

\subsection{Geometric interpretation of the side-splitting map}
The side-splitting map $\phi_*:R_{\br} \rightarrow R_{\br'}$
defined in the previous section induces a map $\phi^*:M_{\br'}
\rightarrow M_{\br}$ on the moduli spaces of polygons.  This map
has a very natural geometric description, which we now give.

Let $x \in M_{\br'}$ be an $n'$-gon.  Then $\phi^*(x)$ is obtained
from $x$ by splitting the $i$th edge of $x$ into $k$ pieces of
lengths $r_{j_1}, \ldots, r_{j_k}$ respectively, where $j_1,
\ldots, j_k$ are the solutions to $\phi(j)=i$.  By ``splitting''
an edge we simply mean regarding a point along the edge as a
vertex.  Thus if $\phi$ is not the identity map then $\phi^*(x)$
will always have 180 degree internal angles.

Geometrically, Proposition~\ref{sidesplitmap} may be rephrased as
follows:

\begin{proposition}
The map $\phi^*:M_{\br'} \rightarrow M_{\br}$ is a closed
immersion of projective varieties.
\end{proposition}

\subsection{Example of the side-splitting map}
Let $\br=(1,1,1,1,1,1,1,1)$, $\br'=(2,2,2,2)$ and
\begin{center}
\begin{tabular}{|c|c|c|c|c|c|c|c|c|} \hline
$i$ & 1 & 2 & 3 & 4 & 5 & 6 & 7 & 8 \\ \hline $\phi(i)$ & 1 & 1 &
2 & 2 & 3 & 3 & 4 & 4 \\ \hline
\end{tabular}.
\end{center}

The semistandard tableau
$$\begin{array}{|c|c|c|c|} \hline
1 & 2 & 4 & 6 \\ \hline 3 & 5 & 7 & 8 \\ \hline
\end{array}$$
is an element of $R_{\br}$; its image under $\phi_*$ is the
semistandard tableau
$$\begin{array}{|c|c|c|c|} \hline
1 & 1 & 2 & 3 \\ \hline 2 & 3 & 4 & 4 \\ \hline
\end{array} \, .$$

On polygons, $\phi^*$ takes equilateral quadrilaterals of side length
2 to equilateral octagons of side length 1.  A typical example is
illustrated below.
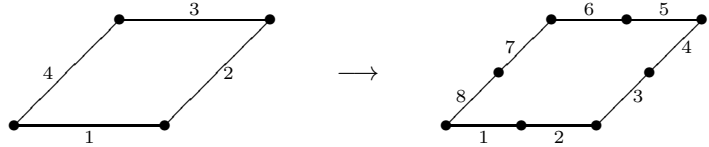
\begin{displaymath}
\begin{xy}
(0,-7)*{}="A"; (20,-7)*{}="B"; (14,7)*{}="C"; (34,7)*{}="D"; "A";
"B"; **\dir{-}; "A"; "C"; **\dir{-}; "B"; "D"; **\dir{-}; "C";
"D"; **\dir{-}; "A"*{\bullet}; "B"*{\bullet}; "C"*{\bullet};
"D"*{\bullet}; (10,-8.5)*{\scriptstyle 1}; (28.5,0)*{\scriptstyle
2}; (24,8.5)*{\scriptstyle 3}; (4.5,0)*{\scriptstyle 4};
\end{xy}
\qquad \longrightarrow \qquad
\begin{xy}
(0,-7)*{}="A"; (20,-7)*{}="B"; (14,7)*{}="C"; (34,7)*{}="D";
(10,-7)*{}="AB"; (7,0)*{}="AC"; (27,0)*{}="BD"; (24,7)*{}="CD";
"A"; "B"; **\dir{-}; "A"; "C"; **\dir{-}; "B"; "D"; **\dir{-};
"C"; "D"; **\dir{-}; "A"*{\bullet}; "B"*{\bullet}; "C"*{\bullet};
"D"*{\bullet}; "AB"*{\bullet}; "AC"*{\bullet}; "BD"*{\bullet};
"CD"*{\bullet}; (5,-8.5)*{\scriptstyle 1}; (15,-8.5)*{\scriptstyle
2}; (25.5,-3)*{\scriptstyle 3}; (32,3.5)*{\scriptstyle 4};
(29,8.5)*{\scriptstyle 5}; (19,8.5)*{\scriptstyle 6};
(8.5,3.5)*{\scriptstyle 7}; (2,-3)*{\scriptstyle 8};
\end{xy}
\end{displaymath}

\subsection{Generators for $R_{\br}$}
The following theorem of Kempe gives generators of the ring
$R_{\br}$ in the case $\br=(1,1,\ldots,1)$ and $n$ is even.  We give a possibly new proof
of this important theorem here. For Kempe's original proof, see
\cite[pg.~156]{Howe}.
We give another proof because we need to extend it to the case of
arbitrary weights (Theorem~\ref{generatorstheorem}).
The graphical description of the ring used here will
be key to results of \cite{HowardMillsonSnowdenVakil}.

\begin{theorem}[Kempe, 1894]
\label{degonegen} For $\br=(1,1,\ldots,1)$ and $n$ even, the ring $R_{\br}$ is
generated by degree one tableaux; in other words, 
$R_{\br}$ is generated by $2 \times (n/2)$ semistandard tableaux
for which each integer between $1$ and $n$ occurs precisely once.
\end{theorem}

\begin{proof}
We associate a graph $\Gamma$ to a tableau $\tau$ as follows:  the vertices of the
graph are the integers between $1$ and $n$; for each column in the
tableau containing the number $i$ and $j$, place an edge between
$i$ and $j$.  The graph is undirected (we have assumed that $i < j$ since 
$\tau$ is semistandard) and allowed to contain
multiple edges between vertices.

A tableau $\tau$ can be written as a product $\tau_1 \tau_2$ of tableaux
where $\tau_1$ is degree 1 if and only if one can pick $n/2$ edges in
the graph associated to $\tau$ which contain each vertex precisely
once.  In fact, these edges can be used to form the columns of
$\tau_1$ by a process inverse to that of the previous paragraph.  If
$\tau$ can be decomposed as such we will say that it (or its graph)
can be \emph{factored}.

We now examine how Pl\"ucker relations work on graphs.  On
tableaux, the Pl\"ucker relation takes a tableau and two columns
in that tableau and yields a sum of two tableaux.  In terms of
graphs, therefore, the Pl\"ucker relation should take a graph and
two edges in the graph and result in two new graphs.  More
precisely, let $\tau$ be a tableau and $\Gamma$ its associated graph.
Let $(a, b)$ and $(c, d)$ be two edges in $\Gamma$ ({i.e.},
two columns in $\tau$).  The Pl\"ucker relations applied to these two
columns in $\tau$ give
$$\begin{array}{|c|c|} \hline a & b \\ \hline c & d \\ \hline \end{array}=
\begin{array}{|c|c|} \hline a & c \\ \hline b & d \\ \hline \end{array}+
\begin{array}{|c|c|} \hline a & b \\ \hline d & c \\ \hline \end{array}.$$
Thus the graphs of the resulting tableaux are obtained from
$\Gamma$ by removing the edges $(a,c)$, $(b,d)$ and then in the
first graph adding edges $(a,b)$,$(c,d)$ and in the second adding
$(a,d)$,$(b,c)$.  This may be represented visually by:
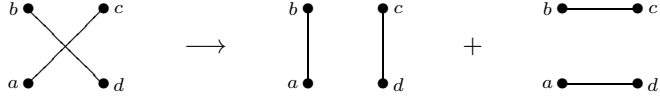
\begin{displaymath}
\begin{xy}
(0,-5)*{}="A"; (0,5)*{}="B"; (10,-5)*{}="D"; (10,5)*{}="C"; "A";
"C"; **\dir{-}; "B"; "D"; **\dir{-}; "A"*{\bullet}; "B"*{\bullet};
"C"*{\bullet}; "D"*{\bullet}; (-2,-5)*{\scriptstyle a};
(-2,5)*{\scriptstyle b}; (12,-5)*{\scriptstyle d};
(12,5)*{\scriptstyle c};
\end{xy}
\qquad \longrightarrow \qquad
\begin{xy}
(0,-5)*{}="A"; (0,5)*{}="B"; (10,-5)*{}="D"; (10,5)*{}="C"; "A";
"B"; **\dir{-}; "C"; "D"; **\dir{-}; "A"*{\bullet}; "B"*{\bullet};
"C"*{\bullet}; "D"*{\bullet}; (-2,-5)*{\scriptstyle a};
(-2,5)*{\scriptstyle b}; (12,-5)*{\scriptstyle d};
(12,5)*{\scriptstyle c};
\end{xy}
\qquad + \qquad
\begin{xy}
(0,-5)*{}="A"; (0,5)*{}="B"; (10,-5)*{}="D"; (10,5)*{}="C"; "A";
"D"; **\dir{-}; "B"; "C"; **\dir{-}; "A"*{\bullet}; "B"*{\bullet};
"C"*{\bullet}; "D"*{\bullet}; (-2,-5)*{\scriptstyle a};
(-2,5)*{\scriptstyle b}; (12,-5)*{\scriptstyle d};
(12,5)*{\scriptstyle c};
\end{xy}
\end{displaymath}

We first show that every
degree two tableaux can be expressed algebraically in terms of
degree one tableaux.  Now, the graph of such a
tableau is a graph on $n$ vertices for which every vertex has
valence 2.  Such a graph is necessarily a disjoint union of
cycles.  If we take two cycles and one edge from each cycle and
apply the Pl\"ucker relations, the cycles merge into one cycle in
both of the resulting graphs.  Thus by applying the Pl\"ucker
relations repeatedly we can merge together pairs of cycles of odd length 
and end up with many graphs each of which consists entirely 
of even length cycles.  An even cycle can clearly
be factored by taking every other edge.  Thus the original degree two 
tableau can be written as a sum of
products of degree one tableaux.

Now suppose that $\tau$ has degree $k > 2$.  Let $\Gamma$ be the associated graph; note 
that $\Gamma$ is $k$-regular. 
We form a new graph $\widetilde{\Gamma}$ by doubling the vertex set; 
each vertex is split into a ``man'' and a ``woman'', i.e., the vertex set of 
$\widetilde{\Gamma}$ is: $$
\t{Vert}(\widetilde{\Gamma}) = \{\t{man}_1,\ldots,\t{man}_n,\t{woman}_1,\ldots,\t{woman}_n\}.$$
For each edge $(i,j)$ of $\Gamma$ let $(\t{man}_i, \t{woman}_j)$ and $(\t{man}_j, \t{woman}_i)$
be edges of $\widetilde{\Gamma}$.  Note that $\widetilde{\Gamma}$ is $k$-regular and bi-partite 
between men and women.  Recall Hall's Marriage Theorem for bipartite graphs: there is a matching 
(a $1$-regular subgraph with the same vertex set) in $\widetilde{\Gamma}$ if and only if any collection of  women  are compatible 
with at least as many men.  This condition is satisfied since $\widetilde{\Gamma}$ is regular.  
Hence there is a matching $\widetilde{\Delta}$ in $\widetilde{\Gamma}$.   Let $\Delta$ be 
the graph on vertices $\{1,\ldots,n\}$ such that for each edge $(\t{man}_i, \t{woman}_j)$ 
in $\widetilde{\Delta}$ we place the edge $(i,j)$ into $\Delta$.  Now $\Delta$ is $2$-regular, and it is almost a subgraph of $\Gamma$.  However, there could be two occurrences of an edge $(i,j)$ in 
$\Delta$ but only one in $\Gamma$; this may happen if both 
$(\t{man}_i,\t{woman}_j)$ and $(\t{man}_j,\t{woman}_i)$ are edges of $\widetilde{\Delta}$.
Nevertheless, this problematic repeated edge in $\Delta$ is a 2-cycle.  Hence if we apply Pl\"ucker relations to rewrite $\Delta$ as a sum of $2$-regular graphs with all even cycles, this $2$-cycle 
$(i,j)$ will be carried along in each term as a common factor.  Finally when we pick out a $1$-regular 
subgraph from each such term, the edge $(i,j)$ will be chosen only once.  Hence we will obtain an 
honest expression for $\Gamma$ as a sum of terms with $1$-regular subgraphs.   Thus each tableau 
of degree $k$ is generated by tableaux of degree one.
\end{proof}

We illustrate the proof given above for $n=6$.  Consider the
tableau
$$\begin{array}{|c|c|c|c|c|c|} \hline
1 & 2 & 3 & 4 & 5 & 6 \\ \hline 2 & 3 & 4 & 5 & 6 & 1 \\ \hline
\end{array}$$
Its graph and a factorization are depicted below
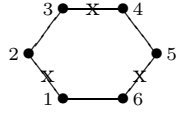
\begin{displaymath}
\begin{xy}
(0,-6)*{}="A"; (-4.5,0)*{}="B"; (0,6)*{}="C"; (8,6)*{}="D";
(12.5,0)="E"; (8,-6)="F"; "A"; "B"; **\dir{-}; "B"; "C";
**\dir{-}; "C"; "D"; **\dir{-}; "D"; "E"; **\dir{-}; "E"; "F";
**\dir{-}; "F"; "A"; **\dir{-}; "A"*{\bullet}; "B"*{\bullet};
"C"*{\bullet}; "D"*{\bullet}; "E"*{\bullet}; "F"*{\bullet};
(-2,-6)*{\scriptstyle 1}; (-6.5,0)*{\scriptstyle 2};
(-2,6)*{\scriptstyle 3}; (10,6)*{\scriptstyle 4};
(14.5,0)*{\scriptstyle 5}; (10,-6)*{\scriptstyle 6}; (-2,-3)*{\rm
x}; (4,6)*{\rm x}; (10.2,-3)*{\rm x};
\end{xy}
\end{displaymath}
This corresponds to the factorization of tableau:
$$\begin{array}{|c|c|c|c|c|c|} \hline
1 & 2 & 3 & 4 & 5 & 6 \\ \hline 2 & 3 & 4 & 5 & 6 & 1 \\ \hline
\end{array}
=
\begin{array}{|c|c|c|} \hline
1 & 3 & 5  \\ \hline 2 & 4 & 6  \\ \hline
\end{array}
\, \cdot \,
\begin{array}{|c|c|c|} \hline
2 & 4 & 6  \\ \hline 3 & 5 & 1  \\ \hline
\end{array}
$$
A more interesting example is given by the tableau
$$\begin{array}{|c|c|c|c|c|c|} \hline
1 & 1 & 2 & 4 & 4 & 5 \\ \hline 2 & 3 & 3 & 5 & 6 & 6 \\ \hline
\end{array}$$
the graph of which is
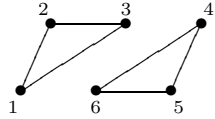
\begin{displaymath}
\begin{xy}
(0,-4)*{}="A"; (4,5)*{}="B"; (14,5)*{}="C"; (24,5)*{}="D";
(20,-4)="E"; (10,-4)="F"; "A"; "B"; **\dir{-}; "B"; "C";
**\dir{-}; "C"; "A"; **\dir{-}; "D"; "E"; **\dir{-}; "E"; "F";
**\dir{-}; "F"; "D"; **\dir{-}; "A"*{\bullet}; "B"*{\bullet};
"C"*{\bullet}; "D"*{\bullet}; "E"*{\bullet}; "F"*{\bullet};
(-1,-6)*{\scriptstyle 1}; (3,7)*{\scriptstyle 2};
(14,7)*{\scriptstyle 3}; (25,7)*{\scriptstyle 4};
(21,-6)*{\scriptstyle 5}; (10,-6)*{\scriptstyle 6};
\end{xy}
\end{displaymath}
Clearly, this cannot be factored.  However, upon applying the
Pl\"ucker relations to edges $(1,3)$ and $(4,6)$ we obtain the two
graphs:
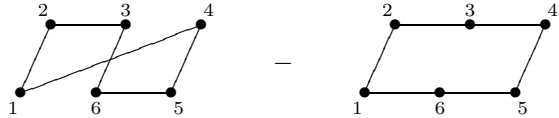
\begin{displaymath}
\begin{xy}
(0,-4)*{}="A"; (4,5)*{}="B"; (14,5)*{}="C"; (24,5)*{}="D";
(20,-4)="E"; (10,-4)="F"; "A"; "B"; **\dir{-}; "B"; "C";
**\dir{-}; "C"; "F"; **\dir{-}; "D"; "E"; **\dir{-}; "E"; "F";
**\dir{-}; "A"; "D"; **\dir{-}; "A"*{\bullet}; "B"*{\bullet};
"C"*{\bullet}; "D"*{\bullet}; "E"*{\bullet}; "F"*{\bullet};
(-1,-6)*{\scriptstyle 1}; (3,7)*{\scriptstyle 2};
(14,7)*{\scriptstyle 3}; (25,7)*{\scriptstyle 4};
(21,-6)*{\scriptstyle 5}; (10,-6)*{\scriptstyle 6};
\end{xy}
\qquad-\qquad
\begin{xy}
(0,-4)*{}="A"; (4,5)*{}="B"; (14,5)*{}="C"; (24,5)*{}="D";
(20,-4)="E"; (10,-4)="F"; "A"; "B"; **\dir{-}; "B"; "C";
**\dir{-}; "C"; "D"; **\dir{-}; "D"; "E"; **\dir{-}; "E"; "F";
**\dir{-}; "F"; "A"; **\dir{-}; "A"*{\bullet}; "B"*{\bullet};
"C"*{\bullet}; "D"*{\bullet}; "E"*{\bullet}; "F"*{\bullet};
(-1,-6)*{\scriptstyle 1}; (3,7)*{\scriptstyle 2};
(14,7)*{\scriptstyle 3}; (25,7)*{\scriptstyle 4};
(21,-6)*{\scriptstyle 5}; (10,-6)*{\scriptstyle 6};
\end{xy}
\end{displaymath}
Each of these graphs is isomorphic to the graph in the previous
example (as unlabelled graphs) and so can be factored.  Such a
factorization would express the original tableau as a sum with two
terms, each of which is a product of two degree one tableaux.

\begin{theorem}\label{generatorstheorem}
For any weight $\br'$ the ring $R_{\br'}$ is generated by lowest
degree elements.
\end{theorem}

\begin{proof}
This follows immediately from  Theorem~\ref{degonegen} and
Proposition~\ref{sidesplitmap}. Take $\br$ to be the
side-length vector of length $n$ consisting of all 1's, where
$n=\vert \br' \vert$.  There is an obvious map $\phi:[1,n]
\rightarrow [1,n']$ satisfying the conditions to be a
side-splitting map.  The resulting surjection $\phi_*:R_{\br}
\rightarrow R_{\br'}$ implies that $R_{\br'}$ is generated by
lowest degree elements, since this is true for $R_{\br}$.
\end{proof}

Note that these results give us projective embeddings of
$(\mathbb{CP}^1)^n\q_\br \SL(2,\C)$: if $N+1$ is the number of
lowest degree elements of $R_{\br}$ then $(\mathbb{CP}^1)^n\q_\br
\SL(2,\C)$ embeds into $\CP^N$.  We will call this embedding the
\emph{Kempe embedding}.  For example, for $n=2m$ even and
$\br=(1,1,\ldots,1)$ the hook length formula gives $N=C_m-1$,
where
$$C_m=\frac{(2m)!}{(m+1)! \, m!}$$
is the $m$th Catalan number; thus the moduli space of equilateral
$2m$-gons naturally embeds nondegenerately into $\CP^{C_m-1}$.


\section{The three isomorphic polytopes}\label{threepolytopes}

In this section, we give three quite different descriptions of the key
polytope used in this paper.  We will use two of them (the diagonal
length polytope and the semistandard tableau polytope) in our
argument, and we discuss the third (the Gel'fand-Tsetlin polytope)
because it is an important description used elsewhere in the
literature.

\subsection{The diagonal length polytope $D(\br)$}

We define the  map $F:\Mr \to \R^{n-3}$ by
$F(\mathbf{e}) = \mathbf{d}$
where $\mathbf{d} = (d_1,d_2,\ldots,d_{n-3})$ is the set of lengths
of the diagonals drawn from the zeroth vertex $v_0$ to the remaining vertices.
Precisely
$$d_i = ||e_1 + e_2 + \cdots + e_{i}||, \ 1 \leq i \leq n-1.$$
Thus the $i$--diagonal joins $v_0$ to $v_{i}$.

We note that $d_1$ and $d_{n-1}$ are fixed (they are the lengths of
the first and last edges), $d_1 = r_1$ and $d_{n-1} = r_n$.  Also $F$
is continuous, but is not smooth where some diagonal length $d_i$ is equal
to zero (because $\sqrt{x}$ is not a smooth function of $x$
near $x=0$, hence $d_i$ is not a smooth function).

We define a  subset $D(\br)$ of $\R^{n-3}$ by
$$D(\br) = F(\Mr).$$

\begin{proposition}
The set $D(\br) = F(\Mr)$ is the convex polytope cut out by the $n-2$
triples of triangle inequalities

\begin{enumerate}
\item $d_i - d_{i+1} \leq r_{i+1}$
\item $d_i - d_{i+1} \geq - r_{i+1}$
\item $d_i + d_{i+1} \geq r_{i+1}$
\end{enumerate}
for $1 \leq i \leq n-2$.

\bigskip

$$
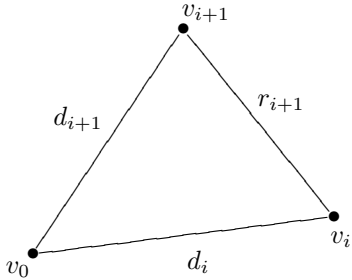
\begin{xy}
(0,0)*{\bullet}="A";
(-2,-2)*{v_0};
(40,5)*{\bullet}= "B";
(41,2)*{v_i};
(20,30)*{\bullet} = "C";
(23,32)*{v_{i+1}};
"A";"B"**\dir{-};
"B";"C"**\dir{-};
"C";"A"**\dir{-};
(22,-1)*{d_i};
(33,20)*{r_{i+1}};
(6,18)*{d_{i+1}};
\end{xy}
$$

\bigskip

\end{proposition}

\begin{proof}
Draw an abstract convex $n$-gon $P$ in the plane (without specifying its
side-lengths). Then the above diagonals
triangulate $P$ into $n-2$ triangles, $T_1, T_2,\ldots,T_{n-2}$ where
$T_1$ has vertices $v_0,v_1$ and $v_2$, $T_2$ has vertices $v_0,v_2,v_3$ and
finally $T_{n-2}$ has vertices $v_0, v_{n-2}$ and $v_{n-1}$.

$$
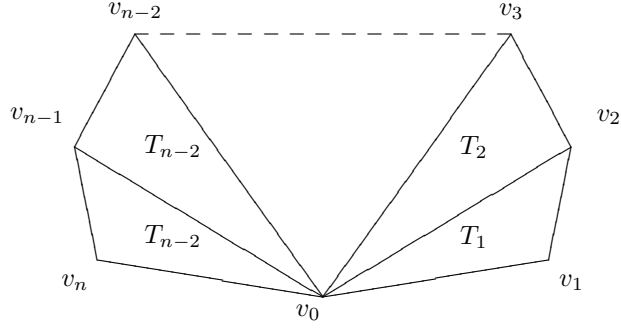
\begin{xy}
(0,0)*{}="A";
(-2,-2)*{v_0};
(30,5)*{}= "B";
(20,8)*{T_1};
(33,2)*{v_1};
(33,20)*{} = "C";
(20,20)*{T_2};
(38,24)*{v_2};
(25,35)*{} = "D";
(25,38)*{v_3};
(-30,5)*{}= "F";
(-20,8)*{T_{n-2}};
(-33,2)*{v_n};
(-33,20)*{} = "G";
(-20,20)*{T_{n-2}};
(-38,24)*{v_{n-1}};
(-25,35)*{} = "H";
(-25,38)*{v_{n-2}};
"A";"F"**\dir{-};
"F";"G"**\dir{-};
"A";"G"**\dir{-};
"G";"H"**\dir{-};
"A";"H"**\dir{-};
"A";"B"**\dir{-};
"B";"C"**\dir{-};
"C";"A"**\dir{-};
"C";"D"**\dir{-};
"D";"A"**\dir{-};
"D";"H"**\dir{--};

\end{xy}
$$

\bigskip

The critical
observation is that we can construct a polygon with the given side-lengths
if and only if we can realize each of the  $n-2$ triangles $T_i, \ 1 \leq
i \leq n-2$ as planar triangles with side-lengths $d_i, r_{i+1}$ and $d_{i+1}$.
For given the $n-2$ triangles we can assemble them along their common diagonals
to obtain the $n$-gon. But the triangle $T_i$ can be realized
if and only if the three numbers $d_i, r_{i+1}$ and $d_{i+1}$ satisfy the
triangle inequalities.
\end{proof}

\subsection{The Gel'fand-Tsetlin polytopes
$GT(\Lambda \varpi_2, \br)$}

To begin this section we will make a definition. (Our presentation
here is borrowed from \cite{DeLoeraMcAllister} though we have indexed
our entries differently.)

For each $n \in \N$, let $X_n$ be the set of all triangular arrays
$(x_{ij})_{1 \le i < j \leq n}$ with $x_{ij} \in \R$.

\begin{definition}
A Gel'fand-Tsetlin pattern or GT-pattern is a triangular array
$\bx = (x_{ij})_{1 \le i < j \leq n} \in X_n$ such that
all $x_{ij}$ are nonnegative and  satisfy the interlacing inequalities:
$x_{i,j+1} \geq x_{ij} \geq x_{i,j+1}$, for $1 \leq i < j \leq n-1$.
\end{definition}

Define $s_i$, $1 \leq i \leq n$, to be the sum of the entries in
the $i$th row,
$$\sum_{j=1}^i x_{ij}$$
and $r_i$, $1 \leq i \leq n$, to be the first differences of the $s_i$'s so
$$r_1 = s_1 \ \text{and} \ r_i = s_{i+1} - s_i, \quad 1 \leq i \leq n-1.$$
We then define

\begin{definition}
The weight $wt(\bx)$ of the Gel'fand-Tsetlin pattern $\bx$ is the $n$-tuple $\br =(r_1,r_2, \ldots,r_n)$ where
the $r_i$'s are as immediately above.
\end{definition}

A GT-pattern is traditionally depicted in an inverted triangle array. See
\cite[pg.~2]{DeLoeraMcAllister}. In our case all the entries $x_{ij}$ will be zero for
$i \geq 3$, that is there are only two nonzero entries in each row
of the pattern (except the first (bottom) row  where there is only one nonzero
entry). We will use $a_i,b_i$ to denote the nonzero entries in the $i$th row.
Furthermore the two nonzero entries in the top row will be assumed
to be equal. Their common value will be denoted $\Lambda$. We give an
example for the case of $n=6$ ( we define $a_n = b_n = a_{n-1}=\Lambda$ and
$b_1 = 0$).

\[
\begin{array}{ccccccccccc}
\Lambda&{}&\Lambda&{}&0&{}&0&{}&0&{}&0\\
{}&\Lambda&{}&b_5&{}&0&{}&0&{}&0&{}\\
{}&{}&a_4&{}&b_4&{}&0&{}&0&{}&{}\\
{}&{}&{}&a_3&{}&b_3&{}&0&{}&{}&{}\\
{}&{}&{}&{}&a_2&{}&b_2&{}&{}&{}&{}\\
{}&{}&{}&{}&{}&a_1&{}&{}&{}&{}&{}
\end{array}
\]
\centerline{Figure 2}

We will often abbreviate such a pattern by the ordered
pair $(\ba,\bb)$ where $\ba = (a_1,a_2,\ldots,a_n)$ and $\bb=(b_1,b_2,\ldots,
b_n)$.

Thus we now have $s_i = a_i + b_i, \ 1 \leq i \leq n$ and
the interlacing inequalities now become

\begin{enumerate}
\item $ a_{i+1} \geq a_i$
 \item $a_i \geq b_{i+1}$
 \item $b_{i+1} \geq b_i.$
\end{enumerate}
Here $1 \leq i \leq n $.

We will denote the the above polytope by $GT(\Lambda \varpi_2)$.
Thus the polytope $GT(\Lambda \varpi_2)$ has affine hull $X(\Lambda)$,
the subspace of $\R^{2n}$ defined by the conditions $a_n =b_n = \Lambda$ and
$b_1 = 0$.

\subsubsection{The polytope $GT(\Lambda \varpi_2, \br)$}

Suppose now that an $n$--tuple of positive real numbers $\br$
is given. We define
$$s_i = \sum_{j=1}^i r_j.$$

We define an affine subspace $X(\Lambda,\br)$ of $\R^{2n}$
by the conditions
$$s_i = a_i + b_i, \qquad 1 \leq i \leq n.$$

\begin{definition} The polytope  $GT(\Lambda \varpi_2, \br)$
is the intersection of the polytope $GT(\Lambda \varpi_2)$
with the affine subspace $X(\Lambda,\br)$.
\end{definition}

Thus the elements $x \in GT(\Lambda \varpi_2, \br)$ are required to
have weight $\br$. Recalling that $b_1= 0$  we have $a_1 = s_1 = r_1$.

\begin{definition}
We define a map
$\Phi:\R^{2n} \to \R^{n-1}$ by $\Phi((\ba,\bb)) = \bd$ where
$$d_i = a_i - b_i$$
for all $i$ such that $1 \leq i \leq n-1$.
\end{definition}

The map $\Phi$ is the map of momentum polytopes underlying
the Hausmann-Knutson duality between the bending integrable system on
the polygon space $M_{\br}$ and the Gel'fand-Tsetlin integrable system on
the torus quotient $\Gr_2(\C^n)\q _{\br}T$, see \cite{HausmannKnutson}.

\begin{proposition}
Suppose $\Lambda = (\sum_{i=1}^n r_i)/2$. Then
the map $\Phi$ carries the Gel'fand-Tsetlin polytope $GT(\Lambda \varpi_2, \br)$
bijectively onto the diagonal length polytope $D(\br)$.
\end{proposition}

\begin{proof}
Note first that the inequality $a_i \geq b_i$ corresponds under $\Phi$
to the inequality $d_i \geq 0$.
We will prove that the interlacing inequalities relating the $i$th and
$i+1$ rows (i.e.\ the rows $(a_i,b_i)$ and $(a_{i+1}, b_{i+1}))$ correspond under
$\Phi$ to the triangle inequalities for
the triangle $T_{i+2}$ (that is the triangle with side-lengths $d_{i+1},
d_{i+2},r_{i+3}$).

The three interlacing inequalities are
\begin{enumerate}
\item $ a_{i+1} \geq a_i$
 \item $a_i \geq b_{i+1}$
 \item $b_{i+1} \geq b_i.$
\end{enumerate}

 The reader will verify that each of the interlacing inequalities corresponds
 to one of the three triangle inequalities. We label the triangle inequalities
 so that corresponding inequalities have the same label.

\begin{enumerate}
\item $d_{i+1} - d_{i+2} \leq r_{i+1}$
 \item $d_{i+1} + d_{i+2} \geq r_{i+1}$
 \item $d_{i+1} - d_{i+2} \geq -r_{i+1}$
 \end{enumerate}

\end{proof}

We leave the analysis of the two boundary triangles to the reader.


\subsection{The polytopes
$SS(\Lambda \varpi_2, \br)$}

\vskip 12pt

The following definitions are fundamental
in combinatorial representation theory.

\begin{definition}
A filling of a Young diagram by the integers between $1$ and $n$ is said
to be semistandard if the columns are strictly increasing and the rows
are weakly increasing.
\end{definition}

A Young diagram with a semistandard filling by the integers between $1$ and $n$
is called a semistandard Young tableau. Given a semistandard Young tableau
$\tau$ we define its weight $wt(\tau)$ to be the sequence of nonnegative integers
$\br = (r_1,r_2,\ldots,r_n)$ where $r_i$ is the number of $i$'s in the tableau.

\bigskip

We will now restrict to the case in which the underlying Young diagram
is a $2$ by $M$ rectangle.

\begin{definition}
Let $\tau$ be a semistandard tableau obtained by a semistandard filling
of a $2$ by $M$ rectangle filled by the integers between $1$ and $n$.
We let $k_{1j}$, resp.\ $k_{2j}$ denote the number of $j$'s in the
first, resp.\ second row. We define the multiweight $\wt(\tau)$ of $\tau$ to be the $2$ by $n$
matrix with integer entries given by
$ \wt(\tau) = \begin{pmatrix} \bk_1 \\
\bk_2
\end{pmatrix} $
where
\[ \begin{pmatrix} \bk_1 \\
\bk_2
\end{pmatrix}
 = \begin{pmatrix} k_{11} & k_{12}  & k_{13} & \cdots & k_{1,n-1} & k_{1n} \\
                 k_{21} & k_{22}  & k_{23} & \cdots & k_{2,n-1} & k_{2n}
\end{pmatrix}.
\]
\end{definition}
For example if 
$\tau$ =
\begin{tabular}{| c | c | c | c| c|c|}
\hline
1 & 1 & 2 & 4 & 4 & 5  \\ \hline
2 & 3 & 3 & 5 & 6 & 6  \\ \hline
\end{tabular}
then
$\wt(\tau)=
  \begin{pmatrix} 2 & 1  & 0 & 2 & 1 & 0 \\
                   0 & 1  & 2 & 0 & 1 & 2
\end{pmatrix}.$ \vskip12pt
Note also that the (ordinary) weight $wt(\tau)=\br$ is given by
$$r_i = k_{1i} + k_{2i}, \qquad 1 \leq i \leq n.$$
The multiplicities $k_{ij}$ satisfy the following ``stairstep inequalities''
(corresponding to the condition that the columns be strictly increasing):
$$\sum_{j=1}^{i} k_{2j} \leq \sum_{j=1}^{i-1} k_{1j}, \quad \quad 2 \leq i \leq n-1.$$
\begin{center}
\begin{tabular}{| c | c | c | c| c|c|}
\hline
$\bullet$  & $\bullet$  & $\bullet$ & \quad & \quad &\quad  \\ \hline
$\bullet$  & $\bullet$  & $\bullet$ & $\bullet$ & \quad & \quad  \\ \hline
\end{tabular}
\end{center}
\centerline{The third stair-step inequality} \vskip 12pt
They also satisfy the three equations
\begin{equation} \label{rowsums}
\sum_{i=1}^n k_{1i} = \sum_{i=1}^n k_{2i}
\end{equation}
\begin{equation}\label{secondrowones}
k_{21} = 0
\end{equation}
\begin{equation}\label{firstrowns}
 k_{1n} = 0.
 \end{equation}

\begin{proposition}
The stairstep inequalities and the equations (\ref{rowsums}), (\ref{secondrowones}), (\ref{firstrowns}), give  conditions
on
$\begin{pmatrix} \bk_1 \\
\bk_2
\end{pmatrix}$
that are necessary and sufficient in order that there exist a semistandard
filling of a $2$ by $n$ rectangle such that there are $k_{1i}$ $i$'s,\ $ 1 \leq i \leq n$,
in the first row and $k_{2i}$ $i$'s, \ $1 \leq i \leq n$, in the second row.
\end{proposition}
\begin{proof}
There is only one way to fill in the rows with the given multiplicities
since the rows are weakly increasing. The $i$th stairstep inequality
gives that the string of $i$'s in the bottom row is completed before the
string of $i$'s in the top row is completed. This guarantees that the
columns are strictly increasing.
\end{proof}

To save space we will replace the column vector notation $\begin{pmatrix} \bk_1 \\
\bk_2
\end{pmatrix}$
by the row vector notation $(\bk_1,\bk_2)$

\bigskip
We now define the third polytope $SS(\Lambda \varpi_2, \br)$ we will need in
this paper.

\begin{definition}
Define the polytope $SS(\Lambda \varpi_2)$ to be the set of $2$ by $n$  matrices with
{\it real} number entries
\[  \mathbf{y}
 = \begin{pmatrix} y_{11} & y_{12}  & y_{13} & \cdots & y_{1,n-1} & y_{1n} \\
                 y_{21} & y_{22}  & y_{23} & \cdots & y_{2,n-1} & y_{2n}

\end{pmatrix},
\]
such that the $y_{ij}$'s
satisfy the stairstep inequalities and the three equalities above. 
We define the polytope $SS(\Lambda \varpi_2, \br)$ to be the subset of
$SS(\Lambda \varpi_2)$ consisting of those $\mathbf{y}$ which are of weight $\br$ i.e.\ such that
$$y_{1j} + y_{2j} = r_j, \ 1 \leq j \leq n.$$
\end{definition}

\subsubsection{The map $\Psi$}

 We define a map $\Psi$ by
 $\Psi(\ba,\bb) = (\bk_1,\bk_2)$ where
 $$k_{1,i+1} = a_{i+1} - a_i \ \text{and} \ k_{2,i+1} = b_{i+1} -b_i.$$

\begin{proposition} \hfill
\begin{enumerate}
\item[(i)] The map $\Psi$ carries the Gel'fand-Tsetlin polytope
$GT(\Lambda \varpi_2))$ bijectively onto the polytope $SS(\Lambda \varpi_2)$.
\item[(ii)] The map $\Psi$ is weight-preserving and consequently carries the Gel'fand-Tsetlin polytope
$GT(\Lambda \varpi_2), \br)$ bijectively onto the polytope
$SS(\Lambda \varpi_2, \br)$ .
\end{enumerate}
\end{proposition}
\begin{proof}
We first prove (i).
We claim that the $i$th stairstep inequality is equivalent under $\Psi$ to the middle interlacing
inequality $a_i \geq b_{i+1}$. Indeed the sum of the first $i$ entries  in the first row of
$\Psi(x) = \begin{pmatrix} \bk_1 \\
\bk_2
\end{pmatrix}$ telescopes to $a_i$ and the sum of the first $i+1$ entries in the second row
telescopes to $b_{i+1}$. The other two interlacing inequalities are equivalent to
the conditions that $k_{1,i+1}$ and $k_{2,i+1}$ are nonnegative.

 \smallskip

To prove (ii) we first claim that under $\Psi$ the number $r_1$ goes to $k_{11}$
and $r_i$ goes to $k_{1,i} + k_{2,i}$ (the number of times $i$ appears in the tableaux
$T$ corresponding to $\Psi(x)$).
Indeed the first statement is obvious and for $i\geq 2$ we have
$r_i = s_i - s_{i+1} = (a_i+b_i) - (a_{i-1}+ b_{i-1}) = (a_i - a_{i-1}) + (b_i +b_{i-1})=
k_{1i} + k_{2i}$.

From this it is immediate that the weight of the pattern $x$ is equal to the weight of
the tableau corresponding to $\Psi(x)$.
\end{proof}

 \begin{remark}
The map on integral points induced by the map $\Psi$ is a special case
of a basic map
in combinatorial representation theory. The integral GT patterns
and the semistandard Young tableau  index weight bases for the
irreducible representation with highest weight $\Lambda \varpi_2$
and the map $\Psi$ indexes a weight-preserving change of basis (see \cite[Fig.~1]{DeLoeraMcAllister}).
\end{remark}

\subsection{Lattice points in the diagonal polytopes $D(\br)$} \hfill
We have defined the polytope isomorphisms
$$\Psi : GT(\Lambda \varpi_2, \br) \to SS(\Lambda \varpi_2, \br),$$
$$\Phi : GT(\Lambda \varpi_2, \br) \to D(\br).$$
Each of these takes integral integral points to integral points.  Furthermore
$\Psi^{-1} : SS(\Lambda \varpi_2, \br) \to GT(\Lambda \varpi_2, \br)$ carries
integral points to integral points.
The integral points of $SS(\Lambda \varpi_2, \br)$ (and hence $GT(\Lambda \varpi_2, \br)$)
are of fundamental importance since they
correspond to semistandard tableaux.

The map $\Phi^{-1}$ does not take integral points to integral points.
There are generally more integral points in $D(\br)$ than
in $GT(\Lambda \varpi_2, \br)$.
In fact, the integral points in $D(\br)$ are in one-to-one correspondence with
the integral points in $GT(2 \Lambda \varpi_2, 2 \br)$.
A simple example is $\br = (1,1,1,1)$.
The Gel'fand-Tsetlin polytope $GT(2 \varpi_2, \br)$ has two integral points corresponding
to the tableaux $\tableau{1}{2} \tableau{3}{4}$ and
$\tableau{1}{3} \tableau{2}{4}$, but $D(\br)$ has three integral points,
$(1,2,1)$, $(1,1,1)$, and $(1,0,1)$.  The point $(1,1,1)$ is not the image of an integral
point in $GT(2 \varpi_2, \br)$.  Since we are primarily concerned with the
semistandard tableaux which generate $R_\br$, which shall say the a point in $D(\br)$ is a
lattice point iff it maps to an integral point under $\Phi^{-1}$.

\begin{proposition}
Denote the lattice points in $D(\br)$ by $D(\br)(\Z)$.
Let $\bd \in D(\br)$.
Then,
$$\bd \in D(\br)(\Z) \iff d_j \equiv \Bigg(\sum_{i=1}^j r_i\Bigg) \!\! \mod 2 \mbox{ for every } j \leq n-1.$$
(In particular if $\br = (1,1,\ldots,1)$ then $\bd \in D(\br)$ is a lattice point iff
$d_i \equiv i \mod 2$ for each $i \leq n-1$.)
\end{proposition}

\begin{proof}
We have that $\bd \in D(\br)(\Z)$ if and only if $\Phi^{-1}(\bd)$ is integral by definition.
This means that $a_j = ((\sum_{i=1}^j r_i) + d_j)/2$ and $b_j = ((\sum_{i=1}^j r_i) - d_j)/2$ are integers
for each $j \leq n-1$.  This is equivalent to $(\sum_{i=1}^j r_i) + d_j \in 2 \Z$, which is equivalent to
$d_j \equiv (\sum_{i=1}^j r_i) \mod 2$.
\end{proof}

\subsection{Application:  $(\mathbb{CP}^1)^n\q_\br \SL(2,\C)$ is rarely a complete
intersection}

\label{application}As an application of this identification of polytopes,
we show that the moduli space is almost never a complete intersection.
We do this by showing that its degree is too low compared to its
codimension.  The main point of this section is to show that
the semistandard tableau polytope  interpretation allows us
to effectively bound the number of points in the polytope, which
yields a previously inaccessible geometric consequence.

Let $H_X$ be the Hilbert polynomial of a projective variety, $X$ and
$L(H_X)$ its leading coefficient. (We will not use this notation
beyond this subsection.)  If $X \subset \CP^N$ be a complete
intersection which is nondegenerate (not contained in any hyperplane),
then the degree of each defining equation of is at least $2$, so the
degree of $X$ is at least $2^{\mathrm{codim} X}$.
As the leading coefficient of the Hilbert polynomial
is the degree divided by $\dim X !$, we have
\begin{equation}
\label{CIeq}
\frac{1}{(\dim{X})!} \le L(H_X) \, 2^{- \mathrm{codim} \, X}.
\end{equation}

We now apply this observation to show that the Kempe embedding is
rarely a complete intersection.  The even and odd cases must be treated
differently.

First, let $n$ be even, $\br=(1,1,\ldots,1)$ and $N=C_{n/2}-1$. The Kempe
embedding gives an embedding of $(\mathbb{CP}^1)^n\q_\br \SL(2,\C)$
into $\CP^N$.  

\begin{theorem}
\label{CIthm} For $n \ge 8$ even, the image of the Kempe embedding
$M_{\br} \rightarrow \CP^N$ is not a complete intersection. For
$n=2,4,6$ the image is a complete intersection.
\end{theorem}

We wish to apply \eqref{CIeq}.  To do so, we need to
know about the leading coefficient of the Hilbert polynomial $H_R$
of the Kempe embedding.

\begin{lemma}
The leading coefficient $L(H_R)$of the Hilbert polynomial of the
ring $R_{\br}$ is $\le 1$.
\end{lemma}

\begin{proof}
The $\ell$th graded component of $R_{\br}$ has for a basis the
semistandard tableaux of weight $\ell$, which correspond bijectively
to integral points in the polytope $\ell SS(\frac{n}{2}\varpi_2,\br)$
(see \S \ref{threepolytopes}), the dilate of the polytope
$SS(\frac{n}{2}\varpi_2,\br)$ by the factor $\ell$.  Thus the
dimension of the $\ell$th graded component of $R_{\br}$ is equal to
the number of integral points in $\ell SS(\frac{n}{2}\varpi_2,\br)$.
For large $\ell$ we have
$$H_R(\ell) = \# \left( \ell SS \left(\frac{n}{2}\varpi_2,\br \right)(\Z) \right).$$
Now there is a map $\pi$ from the polytope $\ell SS(\frac{n}{2}\varpi_2,\br)$
to the $n-3$ dimensional cube with side lengths
equal to $\ell$ given by
$$\pi(\mathbf{k}) = (k_{1,2}, k_{1,3},\cdots,k_{1,n-2}).$$
The reader will verify that $\pi$ is an injection whence
$$H_R(\ell) = \# \left( \ell SS \left( \frac{n}{2}\varpi_2,\br \right)(\Z) \right) \leq (\ell +1)^{n-3}.$$
But the degree of $H_R$ is $n-3$. The previous inequality (valid for large
$\ell$) implies that the leading coefficient of $H_R$ is less than or equal
to $1$.
\end{proof}

We now prove the Theorem.

\begin{proof}[Proof (of Theorem \ref{CIthm})]
Assume the image of the Kempe embedding is a complete
intersection. We now know the following:  $L(H_{M_{\br}}) \le 1$,
$\dim{M_{\br}}=n-3$ and $\mathrm{codim} \, {M_{\br}} =
C_{n/2}-n+2$.  The inequality \eqref{CIeq} would
then give $2^{C_{n/2}-n+2} \le (n-3)!$; however, this is violated
for $n  \ge 8$:  take $\log_2$ of both sides, together with
$\log_2((n-3)!) < (n-3)(n-2)/2$ for $n \geq 7$.  Also
$C_{n/2} - n + 2 > (n-3)(n-2)/2$ for $n \geq 10$.  When $n=8$ we have
$2^{C_4 - 8 + 2} = 2^8 = 256 > 120 = (8-3)!$.

For $n = 2$ the ring $R_\br \cong \C[x]$ and the space (which is a
point) is a complete intersection. For $n=4$ the space
$(\mathbb{CP}^1)^n\q_\br \SL(2,\C)$ is simply $\CP^1$ (and the
Kempe embedding is surjective onto $\CP^1$). For $n=6$ the image of the Kempe
embedding is a cubic hypersurface in $\CP^4$, see \cite{DO}. These
are all complete intersections.
\end{proof}

Now let us examine the case $n$ is odd.  Now the image of the
Kempe embedding lies within $\CP^{N-1}$ where $N$ is the number of
semistandard tableaux weighted by $(2,2,\ldots,2) \in \Z^n$.

\begin{theorem}\label{CIthm_odd}
Let $n \geq 1$ be an odd integer.  Then the image of the Kempe
embedding $M_{\br} \rightarrow \CP^N$ is a complete intersection
iff $n = 1$ or $n = 3$.
\end{theorem}

\begin{lemma}
Let $\br = (2,2,...,2) \in \Z^n$ where $n$ is odd.
The leading coefficient of the
Hilbert polynomial of the ring $R_{\br}$ is $\le 2^{n-3}$.
\end{lemma}

\begin{proof}

The $\ell$th graded component of $R_{\br}$ has for a basis the
semistandard tableaux of weight $2\ell$, which correspond bijectively
to integral points in the polytope $\ell SS(n\varpi_2,\br)$, the
dilation of the polytope $SS(n\varpi_2,\br)$ by
a factor of $\ell$.  Thus the dimension of the $\ell$th graded component
of $R_{\br}$ is equal to the number of integral points in
$\ell SS(n\varpi_2,\br)$.
We repeat the previous argument replacing the cube of dimension $n-3$
with side-lengths $\ell$ by the cube of dimension $n-3$ with
side-lengths $2\ell$.
We obtain
$$H_R(\ell) = \#(\ell SS(n\varpi_2,\br)(\Z)) \leq (2\ell +1)^{n-3}.$$

\end{proof}
We now prove the theorem for the case $n$ odd.

\begin{proof}[Proof (of Theorem \ref{CIthm_odd})]
Assume the image of the Kempe embedding is a complete
intersection. Let $R(n) = N$ be the number of integral points in
$SS(n \varpi_2, \br)$ (the notation $R(n)$ is used because $R(n)$
is the $n$th Riordan number). We now know the following:
$L(H_{M_{\br}}) \le 2^{n-3}$, $\dim{M_{\br}} = n-3$ and
$\mathrm{codim} \, {M_{\br}} = R(n)-n+2$. The inequality 
\eqref{CIeq} would then give $2^{R(n)-n+2} \le
2^{n-3}(n-3)!$.  We claim this inequality is violated for $n \geq
7$.  It is simple to check that $R(n) > 2^{n-3}$ for $n \ge 7$.
But for $n \geq 7$,
$$2^{R(n)-n+2} > 2^{2^{n-3}-n+2} > 2^{n-3} (n-3)!.$$

We then deal with the small cases by hand.
If $n = 5$ then $R(5) = 6$ and $2^{R(n) - n + 2} = 8$ and $2^{n-3}
(n-3)! = 8$ so the inequality is not violated. However, the ring
$R_\br$ is still not a complete intersection; there are six
generators and five essential quadratic relations (this
is classical, and may also be seen in \S
\ref{toricring}). But the space is two dimensional, so to be a
complete intersection its ideal should have just three generators.
Therefore $R_\br$ is not a complete intersection for $n \geq 5$.

If $n = 1$ then $R_\br \cong \C$ (the space is empty) and if $n =
3$ then $R_\br \cong \C[x]$ (the space is a point). These are both
complete intersections.
\end{proof}


\section{A toric degeneration of $M_\br$}\label{toric_degeneration}

We construct a toric degeneration of the moduli space of polygons
by descending a toric degeneration of the Grassmannian
$\Gr_2(\C^n)$ to its torus quotients.
The degeneration of the Grassmannian is essentially the same
as that
given in \cite{LakshmibaiGonciulea} and our toric fiber coincides with that
of  \cite{Sturmfels}, page 104.  Foth and Hu \cite{FothHu}
first observed that the toric degenerations of flag varieties
constructed by Alexeev and Brion \cite{AlexeevBrion} descend to give
toric degenerations of their torus quotients.
The momentum
polytopes of the corresponding toric varieties will be
$SS(\Lambda \varpi_2)$ and $SS(\Lambda \varpi_2,\br)$
respectively.  In what follows we consider problems concerning
lattice points in these polytopes. In all cases the underlying
lattice will be the standard integer lattice.

\subsection{A toric degeneration of $\Gr_2(\C^n)$}

Our toric degeneration is essentially the same as that given in
\cite{LakshmibaiGonciulea}. Recall that $2a = \sum_i r_i = |\br|$ and
$\mathcal{L}_{2,n}^a$ is the very ample line bundle of $\Gr_2(\C^n)$
corresponding to the character $\det^a$ of $\GL(2,\C)$.
Let
$$R =\bigoplus_{N=0}^\infty \Gamma(\Gr_2(\C^n),
(\mathcal{L}_{2,n}^a)^{\otimes N}).$$
The degree one ($N=1$)
sections $s_\tau$ generate $R$ as a ring, where $\tau$ is a semistandard
$2$ by $\Lambda$ tableau filled with indices $1$ through $n$.
For ease of notation, we shall identify $s_\tau$ with $\tau$.

We have that $R$ is an infinite-dimensional $\C$-vector space, with
basis consisting of semistandard $2$ by $Na$ tableaux as $N$ ranges from
$0$ to infinity. (When $N=0$ the section $s_\emptyset$ of the empty tableau is taken to be the
constant section $1$ of the trivial line bundle.)
The multi-weights
$\wt(\tau)$ defined in the previous section make $R$ into a
multi-graded vector space. For each multi-weight $\bk$, let
$\pi_\bk$ denote the projection of $R$ onto the one-dimensional
subspace of weight $\bk$. Let $\gr_1(R)$ denote the associated
multi-graded ring, defined as follows. As a $\C$-vector space $\gr_1(R)$ is the
same as $R$; only the ring structures will differ. We define
multiplication on $\gr_1(R)$ by
$$\bar{\tau} \cdot_{\gr_1(R)} \bar{\sigma} = \pi_\bk(\tau \cdot_R \sigma),$$
where $\bk = \wt(\tau) + \wt(\sigma)$ and $\tau$,$\sigma$ are
semistandard tableaux. Here we have used $\bar{\tau}$ and $\bar{\sigma}$
to denote the elements in $\gr_1(R)$ corresponding to the (multihomogeneous)
elements $\tau$ and $\sigma$.

\begin{definition}
Given semistandard tableaux $\tau$ and $\sigma$, let
$\tau * \sigma$ be the unique semistandard tableau such that
$\wt(\tau * \sigma) = \wt(\tau) + \wt(\sigma)$.
\end{definition}

\begin{remark}
The tableau $\tau * \sigma$ can be described alternatively as the
tableau obtained by concatenating $\tau$ and $\sigma$, then
rearranging the top row and bottom row indices so that both rows are
nondecreasing.
\end{remark}

\begin{lemma}
The product of $\bar{\tau}$ and $\bar{\sigma}$ in $\gr_1(R)$ is given by:
$$\tau \cdot_{\gr_1(R)} \sigma = \overline{\tau * \sigma}.$$
\end{lemma}

\begin{proof}
If basic
straightening relations are applied to the product $\tau \sigma$,
all the monomials have the same set of indices as $\tau$ and $\sigma$
together.
Furthermore, there is exactly one monomial in the sum which has
the same set of
indices in its top row as the set of indices in the top row of the
concatenated $\tau \sigma$. Finally when enough applications of
basic straightening relations are applied so that each term in the
sum is semistandard, the unique term with the same indices in the top
row as $\tau \sigma$ must be equal to $\tau * \sigma$.
\end{proof}

Let $\widetilde\T$ be the complex torus $\C^{2 \times n}/\Z^{2
\times n}$ of dimension $2n$. Identify the characters of
$\widetilde\T$ with the multi-weights $(\bk_1, \bk_2) \in \Z^{2n}$
by $$(\bk_1,
\bk_2)((t_{1,1},\ldots,t_{1,n},t_{2,1},\ldots,t_{2,n})) =
\prod_{i,j} t_{i,j}^{k_{i,j}}.$$  We define an action of
$\widetilde\T$ on semistandard tableaux. If $\wt(\tau) = (\bk_1,
\bk_2)$ and $t \in \widetilde\T$, let $$t \cdot \tau =
\wt(\tau)(t) \, \tau.$$  Extend the action to be linear on $R$.

The torus $\widetilde\T$ does not act effectively on $R$.
For $k=1,2$ and $1 \leq \ell \leq n$
let $T_{k,\ell}$ denote the one dimensional subgroup of $\widetilde\T$
where all components $t_{i,j}$ except $t_{k,\ell}$ are equal to $1$.
The subgroups  $T_{2,1}$ and $T_{1,n}$ act trivially since in any
semistandard tableau the index $1$ does not appear in the second
row nor does the index $n$ appear in the first row. Furthermore,
since the tableaux also have the same number of entries in the
first and second rows, the one-dimensional subgroup
$(t,t,\ldots,t,t^{-1},t^{-1},\ldots,t^{-1})$ acts trivially.   Let
$\T$ denote the quotient of $\widetilde\T$ by these three
one-dimensional subgroups.  Now the character lattice
$\chi^\ast(\T)$ is canonically a sub-lattice of
$\chi^\ast(\widetilde\T)$ of dimension $2n-3$, given by the
equations
$$k_{2,1} = k_{1,n} = 0 \quad \text{and} \quad \sum_\ell k_{1,\ell} =
\sum_\ell k_{2,\ell}.$$  Now $\T$ acts effectively on $R$ (and
$\gr_1(R)$) as a collection $\C$ module homomorphisms.
Furthermore, $\T$ acts by ring homomorphisms on $\gr_1(R)$, since
$$(t \cdot \tau)(t \cdot \sigma) = \wt(\tau)(t) \wt(\sigma)(t)
\tau \sigma = \wt(\tau * \sigma)(t) (\tau*\sigma) = t \cdot (\tau
* \sigma).$$
Let $\Gr_2(\C^n)_0$ denote the associated toric variety
$\t{Proj}(\gr_1(R))$, which contains the projective quotient $\T /
\Delta(\T)$ of $\T$ as an open subset. Here $\Delta(\T)$ denotes
the image in $\T$ of the diagonal subgroup of $\widetilde\T$.

Now we shall construct $\Gr_2(\C^n)_0$ as the special
fiber of a degeneration of the ring $R$ by choosing a $\T$ stable
$\N$-filtration as in \cite{AlexeevBrion}.
\begin{definition}\label{LGgrading}
Let $C \geq 2$ be an integer. For each semistandard $2$ by $Na$
tableau $\tau = \tableau{\tau_{1,1}}{\tau_{2,1}}
\tableau{\tau_{1,2}}{\tau_{2,2}} \cdots
\tableau{\tau_{1,Na}}{\tau_{2,Na}}$ define
$$C_\tau = \sum_{\ell = 1}^m  \tau_{1,\ell} + C \tau_{2,\ell}.$$
For each non-negative integer $m$ let
$$F_m(R) = \bigoplus_{C_\tau \leq m} \! \! \C[\tau] \subset
R.$$ Let $\gr_2(R)$ be the associated graded ring, with graded
components $F_m(R) / F_{m-1}(R)$.  We shall say that the tableau
$\tau$ has LG-degree $C_\tau$ (since this is the construction of
Lakshmibai-Gonciulea).
\end{definition}

Now $R$ is both filtered \emph{and} graded.  The filtration is given by
the LG-degrees of tableaux and the grading comes from the
decomposition of $R$ into sections of tensor powers of
$\mathcal{L}_{2,n}^a$.  The filtration and grading are not compatible,
since the set of elements of $R$ of degree $k$ and less is not a union
of $LG$--filtration levels.

Now $\gr_2(R)$ is graded in two ways;
we have the standard degree $\deg(\tau) = N$ where $\tau$ is
$2$ by $Na$ and there is also the LG-degree $C_\tau$.  The reason why the
standard grading (by the standard degree) is well-defined on $\gr_2(R)$ is that
the standard degree of a tableau may be computed from its multi-weight.

We will use the following notation for the rest of the paper. We have
an identification of vector spaces
$$ R = \bigoplus \C \tau_i$$
where $\tau_i$ runs over the set of semistandard tableau with two rows.
If the LG-degree of $\tau_i$ is $m$ then we let $\bar{\tau}_i$
denote the corresponding (basis) element of $\gr_2(R)^{(m)} = F_m(R)/F_{m-1}(R)$.
Thus as vector spaces we have
$$\gr_2(R) = \bigoplus \C \bar{\tau}_i.$$

\begin{lemma}
Multiplication within $\gr_2(R)$ is given by
$\bar{\tau} \cdot_{\gr_2(R)} \bar{\sigma} = \overline{\tau * \sigma}$.
\end{lemma}

\begin{proof}
Suppose that $\tau$,$\sigma$ are semistandard tableaux.
Let $\tau \sigma = \tau * \sigma + \sum_i c_i \gamma_i$
be the product of $\tau$
and $\sigma$ computed in $R$, where the $\gamma_i$ are semistandard.
We shall show that $C_{\gamma_i} < C_{\tau * \sigma}$ for each $i$.
Consider the basic straightening relations
$${\tableau{i_1}{i_4}} \cdot {\tableau{i_2}{i_3}}
= {\tableau{i_1}{i_3}} \cdot {\tableau{i_2}{i_4}} -
{\tableau{i_1}{i_2}} \cdot {\tableau{i_3}{i_4}}$$
where $1 \leq i_1 < i_2 < i_3 < i_4 \leq n$. (We are considering the above
$2$ by $2$ tableaux as factors of larger tableaux which actually lie within $R$.) Let
$$\alpha = {\tableau{i_1}{i_4}} \cdot {\tableau{i_2}{i_3}}\, , \quad
\alpha_1 = {\tableau{i_1}{i_3}} \cdot {\tableau{i_2}{i_4}}\, ,
\quad \alpha_2 = {\tableau{i_1}{i_2}} \cdot {\tableau{i_3}{i_4}}
\, .$$ We have
$$C_{\alpha} = i_1 + i_2 + C(i_3 + i_4),$$
$$C_{\alpha_1} = i_1 + i_2 + C(i_3 + i_4),$$
$$C_{\alpha_2} = i_1 + i_3 + C(i_2 + i_4).$$
Now $C_\alpha = C_{\alpha_1}$. Since $i_3 > i_2$ and $C \geq 2$,
$$C_\alpha - C_{\alpha_2} = (i_2 - i_3) + C(i_3 - i_2) = (C-1)(i_3-i_2)
> 0.$$
Since the equation $\tau \sigma = \tau * \sigma + \sum_i c_i \gamma_i$
is gotten from a sequence of basic straightening relations, and
$\tau * \sigma$ is the final leftmost term, it is clear from the above that
each $C_{\gamma_i} < C_{\tau * \sigma}$.  Hence the product of
$\bar{\tau}$ and $\bar{\sigma}$ in $\gr_2(R)$ is $\overline{\tau * \sigma}$.
\end{proof}

\begin{corollary}
The $\C$--algebras $\gr_1(R)$ and $\gr_2(R)$ are isomorphic.
\end{corollary}

\begin{remark}
Of course $\gr_1(R)$ and $\gr_2(R)$ are not isomorphic
as graded $\C$--algebras: the grading of $\gr_1(R)$ is by a cone in the character
lattice of $\T$, whereas the grading of $\gr_2(R)$ is a grading
by the natural numbers.
\end{remark}

It is well-known that if $R$ is a filtered ring then there is a one-parameter flat degeneration
with special fiber the associated graded ring of $R$.
We sketch one way to do this, borrowed from \cite{AlexeevBrion}, using the Rees algebra.
Let $z$ be an
indeterminate and let $\mathcal{R}$ be the Rees algebra
$$\mathcal{R} = \bigoplus_{m=0}^\infty F_m(R) z^m \subset
R[z].$$

\begin{theorem}(Alexeev--Brion \cite{AlexeevBrion})
\begin{itemize}
\item $\mathcal{R}$ is flat over $\C[z]$. \item $\mathcal{R}
\otimes_{\C[z]} \C[z,z^{-1}] \cong R[z,z^{-1}]$. \item
$\mathcal{R} \otimes_{\C[z]} \C[z]/ (z) \cong \gr_2(R)$.
\end{itemize}
\end{theorem}

\subsection{Restriction of the degeneration to torus invariants}
The small torus $T \cong (\C^\times)^n$ acts on $R$ by ring homomorphisms,
by $$t \cdot \tau = \chi_\br^a(t) \chi_\bs(t^{-1}) \tau,$$
where $\tau$ is a $2$ by $Na$ semistandard tableau and $\bs = (s_1,s_2,\ldots,s_n)$ where
$s_i$ is the number of occurrences of $i$ in $\tau$.
We extend the action of the small torus $T$ to $\mathcal{R}$ by making
the indeterminate $z$ an invariant.  Let $\mathcal{R}^T$ denote
the torus invariants.

\begin{theorem} \hfill
\begin{itemize}
\item $\mathcal{R}^T$ is flat over $\C[z]$. \item $\mathcal{R}^T
\otimes_{\C[z]} \C[z,z^{-1}] \cong R^T[z,z^{-1}]$. \item
$\mathcal{R}^T \otimes_{\C[z]} \C[z]/ (z) \cong \gr_2(R)^T$.
\end{itemize}
\end{theorem}

\begin{proof}
Flatness follows from the fact that restriction to $T$-invariants is an exact functor of
$\C[z]$-modules.  The other properties are immediate.
\end{proof}

Recall that the T-invariant sub-module $R_\br = R^T$ of $R$ is
generated by symbols $\bar{\tau}$ for those tableaux $\tau$,
with weight a
multiple of $\br$.  We still say that the LG-degree of $\bar{\tau}$ in
$R^T$ is $C_\tau$ as above.  However, the standard degrees of $\tau$
and $\bar{\tau}$ are no
longer the number of columns of $\tau$, but rather the number
$N$ whenever $\tau$ has shape $2$ by $Na$ (where $a = |\br|/2$),
since then $\tau$ corresponds to a section $s_\tau$ of the $N$th tensor
power of the line bundle $\mathcal{L}^{\otimes a}_{2,n}$ over
$\Gr_2(\C^n)$, see \S \ref{AlgGeomSpace}.

\begin{definition}
We  let  $(R_{\br})_0$ denote the projective coordinate ring of the
toric fiber of the above induced toric degeneration or $R_\br = R^T$.
Thus
$$ (R_{\br})_0 \cong \gr_2(R)^T.$$
\end{definition}

The main point of what follows is that the graded $\C$--algebra $(R_{\br})_0$
can be identified with the
semigroup algebra of the graded semigroup $S_{\br}$ of lattice points in the family
of integral dilates $D(N\br),N \geq 0$, of the diagonal polytope $D(\br)$.
Here we  define the degree  of a lattice point $\bd \in D(N\br)$ to be
$N$. This definition of degree makes
the semi-group algebra $\C[S_\br]$ into a graded $\C$-algebra with
relations generated by binomial relations $x_{\bd_1} x_{\bd_2}
- x_{\bd_1'} x_{\bd_2'}$ where $\bd_1 + \bd_2 = \bd_1' + \bd_2'$
in $S_\br$. Here we use the symbol $x_{\bd}$ to be the element of
$\C[S_\br]$ corresponding to $\bd \in D(N\br)$.

\begin{proposition}
We have an isomorphism of graded $\C$--algebras
$$(R_\br)_0 \cong \C[S_\br].$$
\end{proposition}

\begin{proof}
First of all $\gr_2(R^T)$ and $R^T$ are isomorphic as $\C$-modules.  A basis
for the $N$th graded summand is given by the set $\bar{\tau}$ for $\tau$
a semistandard $2$ by $Na$ tableaux of weight  $N \br$.   These correspond to
the integral points occurring in $SS(Na \varpi_2, N \br)$ for $N \geq 0$.  The map
$\Psi \circ \Phi^{-1}$ is a degree--preserving
bijection between $S_\br$ and this basis.  We claim this map induces 
a graded ring homomorphism.
Multiplication in $\gr_2(R^T)$ is given by $\bar{\tau} \bar{\sigma} =
\overline{\tau \ast \sigma}$.  But the multi-weight of
$\tau \ast \sigma$ is the sum of the multi-weights of $\tau$ and $\sigma$ by definition.
Hence multiplication in $\gr_2(R^T)$
corresponds to addition of multi-weights, which is the (commutative) semigroup
structure of the integral points in the various $SS(Na \varpi_2, N \br)$.   The map
$\Psi \circ \Phi^{-1}$ is additive and so it preserves the semigroup structure.  Therefore it induces
to a ring homomorphism.
\end{proof}

We will identify the rings $(R_\br)_0$ and $\C[S_\br]$ using the above
isomorphism henceforth.

\section{The projective coordinate ring of the toric fiber
$(M_\br)_0$}\label{toricring}

Fix $\br \in (\Z^+)^n$. We shall show in this section that the
degree one and degree two elements of the toric ring $\C[S_\br]$
generate $\C[S_\br]$. Furthermore we will show that the
ideal of $\C[S_\br]$ (and hence $(R_\br)_0$) is generated by relations of degrees two,
three, and four.  It follows immediately from the proofs given
that $\C[S_{2\br}] \cong (R_{2 \br})_0$ has a presentation by degree one generators
and quadratic relations.

\subsection{Generators for $\C[S_\br]$.}

\begin{definition}
If $A$, $B$ are subsets of a vector space, their
Minkowski sum $A + B$ is
$$A + B = \{a + b \mid a \in A, \, b \in B\}.$$
Additionally when the vector space contains a given lattice let
$A(\Z)$ denote the lattice points which are contained in $A$.
\end{definition}

\begin{definition}
For each positive integer $k$, let $D(k) = D(k \br)$.
\end{definition}

\begin{lemma}\label{odd_degree}
For each positive integer $m$,
$$D(2m+1)(\Z) = D(1)(\Z) + D(2m)(\Z).$$
\end{lemma}

\begin{proof}
Given $\bd \in D(2m+1)(\Z)$ we shall construct $\bd' \in D(1)(\Z)$
in the proximity of $\bd/(2m+1)$ such that $\bd - \bd' = \bd'' \in
D(2m)(\Z)$.  An element $\bd' = (d'_1,\ldots,d'_{n-1})$ in $D(1)$
is a lattice point iff each $d'_i$ is an integer with the same
parity as $(r_1 + \cdots + r_i)$; that is
$d'_i \equiv (r_1 + \cdots + r_i) \, \mod \, 2$.  On the other
hand $\bd'' = (d''_1,\ldots,d''_{n-1})$ in $D(2m)$ is a lattice
point iff each $d''_i$ is an even integer.  Define
$$d'_i = k, \; \mbox{such that} \; k \in \Z,\; k \equiv (r_1 + \cdots + r_i) \mod
2,\; \mbox{and }\Big|k - \frac{d_i}{2m+1}\Big| < 1.$$ For this to
be well-defined we need to show that $k$ exists and is unique.
Uniqueness follows immediately since there can be only one integer
of a given parity in an open interval of length $2$. For existence
we must check that $d_i/(2m+1)$ does not have opposite parity to
$(r_1 + \cdots + r_i)$, since in this case there is no integer of
the correct parity less than one unit from $d_i/(2m+1)$. But $d_i
\equiv (r_1 + \cdots + r_i) \mod 2$ since $\bd \in D(2m+1)(\Z)$ is
a lattice point and $(2m+1)$ is odd. Hence if $d_i/(2m+1)$ is an
integer then $d_i/(2m+1) \equiv (r_1 + \cdots + r_i) \mod 2$ as
well since $2m+1$ is odd.  Therefore the parity condition for
$\bd'$ is satisfied.

Let $\bd'' = \bd - \bd'$.  Note
that each $d''_i$ is even since $d_i$ and $d'_i$ have the same
parity.  We have that in fact $d''_i$ is the nearest even integer
to $2md_i/(2m+1)$.
Thus the lattice point conditions are satisfied.  Now we
only need to show that $\bd' \in D(1)$ and $\bd'' \in D(2m)$.

Since $d_1 = (2m+1)r_1$ and $d_{n-1} = (2m+1)r_{n-1}$ we have that $d'_1 = r_1$,
$d''_1 = 2m r_1$, $d'_{n-1} = r_{n-1}$, and $d''_{n-1} = 2m r_{n-1}$.  It remains to
show the triangle inequalities,
\begin{itemize}
\item[(1)] $d'_{i-1} \leq d'_i + r_i$, $d''_{i-1} \leq d''_i + 2m r_i,$
\item[(2)] $d'_i \leq d'_{i-1} + r_i$, $d''_i \leq d''_{i-1} + 2m r_i,$
\item[(3)] $r_i \leq d'_{i-1} + d'_i$, $2m r_i \leq d''_{i-1} + d''_i.$
\end{itemize}
Let $\ell = 2m+1$. We have that $|d_{i-1} - d_i| \leq \ell r_i$, hence
$|d_{i-1}/\ell - d_i/\ell| \leq r_i$.  Recall that $d'_{i-1}$ is
the nearest integer to $d_{i-1}/\ell$ with parity
$(r_1 + \cdots + r_{i-1}) \mod 2$ and
$d'_i$ is the nearest integer to $d_i/\ell$ with parity $(r_1 + \cdots + r_i) \mod
2$.  The distance between $d_{i-1}/\ell$ and $d_i/\ell$ is
at most $r_i$ and also $d'_{i-1} - d'_i \equiv r_i \mod 2$.
We also have that $|d'_{i-1} - d_{i-1}/\ell| < 1$ and $|d'_{i} - d_{i}/\ell| < 1$.
Therefore $|d'_{i-1} - d'_i| < r_i + 2$ and consequently
$|d'_{i-1} - d'_i| \leq r_i + 1$.  But $|d'_{i-1} - d'_i| \neq r_i + 1$ because
$d'_{i-1} - d'_i \equiv r_i \mod 2$.
It follows that $|d'_{i-1} - d'_i| \leq r_i$.
Also $|d''_{i-1} - \frac{\ell - 1}{\ell}d_{i-1}| < 1$ and
$|d''_{i} - \frac{\ell - 1}{\ell}d_{i}| < 1$.  Therefore
$|d''_{i-1} - d''_i| < (\ell-1)r_i + 2$, so
$|d''_{i-1} - d''_i| \leq (\ell-1)r_i + 1$.  But $d''_{i-1} - d''_i$ is even
and $(\ell - 1)r_i + 1$ is odd since $\ell$ is odd.  Therefore
$|d''_{i-1} - d''_i| \leq (\ell-1)r_i = 2m r_i$.
Therefore both (1) and (2) hold.

We have $d'_{i-1} + d'_i > d_{i-1}/\ell + d_i/\ell - 2 \geq r_i - 2$.
Thus $d'_{i-1} + d'_i \geq r_i - 1$.  But $d'_{i-1} + d'_i \equiv r_i \mod 2$
so $d'_{i-1} + d'_i \neq r_i - 1$.  Thus
$d'_{i-1} + d'_i \geq r_i$.
We have $d''_{i-1} + d''_i > \frac{\ell-1}{\ell}(d_{i-1} + d_i) - 2 \geq (\ell - 1)r_i - 2$.
Thus $d''_{i-1} + d''_i \geq (\ell - 1)r_i - 1$.  But $d''_{i-1} + d''_i$ is even and
$\ell$ is odd so $d''_{i-1} + d''_i \neq (\ell-1)r_i - 1$.  Thus
$d'_{i-1} + d'_i \geq (\ell - 1)r_i = 2m r_i$.
Therefore (3) holds.
\end{proof}

\begin{lemma}\label{even_degree}
For each positive integer $m$,
$$D(2m)(\Z) = \underbrace{D(2)(\Z) + \cdots + D(2)(\Z)}_m.$$
\end{lemma}

\begin{proof}
We show the lemma by induction on $m \geq 1$.  The case $m=1$ is a
tautology.  Suppose that $m \geq 2$.  We show that
$$D(2m)(\Z) = D(2)(\Z) + D(2m-2)(\Z).$$
Suppose $\bd = (d_1,\ldots,d_{n-1}) \in D(2m)(\Z)$.  We construct
$$\bd ' \in D(2)(\Z), \quad \bd ''~\in~ D(2m-2)(\Z)$$ where
$\bd ' + \bd '' = \bd$, by placing $\bd '$ in the proximity of
$\bd/m$.  Recall that the integrality condition is that the
components of $\bd$, $\bd'$, and $\bd''$ are even integers.

Let $e^- : \R \to 2\Z$ be the function which assigns the nearest
even integer, where odd integers $2t+1$ are mapped to $2t$. To be
concise,
$$e^-(x) = \min \{k \in 2\Z : k + 1 \geq x\}.$$
Similarly let $e^+ : \R \to 2\Z$ assign the nearest even integer
where odd integers $2t+1$ are mapped to $2t+2$,
$$e^+(x) = \max \{k \in 2\Z : k - 1 \leq x\}.$$
We will often use the following properties of $e^-$ and $e^+$:
\begin{itemize}
\item each of $e^-$ and $e^+$ is weakly increasing. \item if $k
\in 2\Z$, then $e^\pm(x + k) = e^\pm(x) + k$.  \item $e^+(-x) = -
e^-(x)$. \item if $x+y \in 2\Z$, then $e^+(x) + e^-(y) = x + y$.
\item if $x+y \geq k \in 2\Z$, then $e^+(x) + e^-(y) \geq k$.
\end{itemize}
Let $\bd = (d_1,\ldots,d_{n-1}) \in D(2m)(\Z)$. Let
$$\mathcal{J}_\bd^0 = \{ i \mid d_{i-1}/m \mbox{ and } d_i/m \mbox{ are odd integers},
d_{i-1} + d_i = 2m r_i, \, 2 \leq i \leq n-1\},$$
$$\mathcal{J}_\bd^1 = \{ i \mid d_{i-1} \leq 2m r_i, \, d_i \leq 2m r_i,
\,2 \leq i \leq n-1\}.$$ Clearly $\mathcal{J}_\bd^0 \subset
\mathcal{J}_\bd^1$.  Let $\mathcal{J}_\bd$ be such that
$\mathcal{J}_\bd^0 \subset \mathcal{J}_\bd \subset
\mathcal{J}_\bd^1$. Let $\{i_1,\ldots,i_s\} = \mathcal{J}_\bd$
such that $i_t < i_{t+1}$ for all t, and set $i_0 = 1$ and
$i_{s+1} = n$.

\noindent Let $\bd ' = (d'_1,\ldots,d'_{n-1}) \in (2\Z)^{n-1}$ be
$$ d'_i=\left \{ \begin{array}{ll}
e^-(d_i/m) & \textrm{for $i_{2t} \leq i < i_{2t+1}$}, 2t \leq s,\\
e^+(d_i/m) & \textrm{for $i_{2t+1} \leq i < i_{2t+2}$}, 2t+1 \leq
s.
\end{array} \right.$$
Let $\bd '' = (d''_1,\ldots,d''_{n-3}) \in (2\Z)^{n-1}$ be
$$ d''_i=\left \{ \begin{array}{ll}
e^+((m-1)d_i/m) & \textrm{for $i_{2t} \leq i < i_{2t+1}$}, 2t \leq s, \\
e^-((m-1)d_i/m) & \textrm{for $i_{2t+1} \leq i < i_{2t+2}$}, 2t+1
\leq s.
\end{array} \right.$$
We will show that $\bd ' + \bd '' = \bd$, $\bd' \in D(2)(\Z)$, and
$\bd'' \in D(2m-2)(\Z)$.  Note that $e^\pm(x) + e^\mp(y) = k$
whenever $x + y = k$ and $k \in 2\Z$. We have that $d_i/m +
(m-1)d_i/m = d_i \in 2\Z$ for all $i$, so
$$d'_i + d''_i = e^\pm(d_i/m) +
e^\mp((m-1)d_i/m) = d_i.$$ Thus $\bd ' + \bd '' = \bd$.

We show that $\bd ' \in D(2)$. The proof that $\bd'' \in D(2m-2)$
is similar. Since $d_1 = 2m r_1$ and $d_{n-1} = 2m r_{n-1}$, we have $d'_1 =
2 r_1$ and $d'_{n-1} = 2 r_{n-1}$. Now suppose $2 \leq i \leq n-1$. We must
show the three triangle inequalities that define $D(2)$:
\begin{itemize}
\item[(1)] $d'_i \leq d'_{i-1} + 2 r_i$, \item[(2)] $d'_{i-1} \leq d'_i
+ 2 r_i$, \item[(3)] $2 r_i \leq d'_{i-1} + d'_i$.
\end{itemize}
Suppose that $i \notin \mathcal{J}_\bd$. We have $d'_{i-1} =
e^\pm(d_{i-1}/m)$ and $d'_{i-1} = e^\pm(d_i/m)$ (the same function
is applied to each). The functions
$e^-$ and $e^+$ are weakly increasing, and since $d_i/m \leq
d_{i-1}/m + 2 r_i$ we have that $$d'_i = e^\pm(d_i/m) \leq
e^\pm(d_{i-1}/m + 2 r_i) = e^\pm(d_{i-1}/m) + 2 r_i = d'_{i-1} + 2 r_i$$ so
(1) holds. Similarly inequality (2) holds. Since $i \notin
\mathcal{J}_\bd^0$ we know that either $d_{i-1} + d_i > 2mr_i$ or
one of $d_{i-1}/m$ or $d_i/m$ is not an odd integer.
Suppose $d_{i-1} + d_i > 2mr_i$.  We have that
$d'_{i-1} + d'_i \geq d_{i-1}/m + d_i/m - 2 > 2r_i - 2$.  But
$d'_{i-1} + d'_i$ is even so $d'_{i-1} + d'_i \geq 2r_i$ and (3) holds.
Suppose that one of $d_{i-1}/m$ or $d_i/m$ is not an odd integer and
$d_{i-1} + d_i = 2mr_i$.
Without loss of generality suppose that $d_{i-1}/m$ is not odd.
Then $e^+(d_{i-1}/m) = e^-(d_{i-1}/m)$. Since the sum
$d_{i-1}/m + d_i/m = 2r_i$ is even, we have that
$e^+(d_{i-1}/m) + e^-(d_i/m) = 2r_i$.  Now
$$d'_{i-1} + d'_i \geq e^-(d_{i-1}/m) + e^-(d_{i}/m) =
e^+(d_{i-1}/m) + e^-(d_i/m) = 2r_i,$$
and again (3) holds.

Suppose that $i \in \mathcal{J}_\bd$.  Hence $i \in
\mathcal{J}_\bd^1$ and so $d_{i-1}/m \leq 2r_i$ and  $d_i/m \leq 2r_i$.
Therefore each of $d'_{i-1} \leq 2r_i$ and $d'_i \leq 2r_i$. We have
$d'_{i-1} = e^\pm(d_{i-1}/m)$ and $d'_i = e^\mp(d_i/m)$. Whenever
two numbers $x,y$ satisfy $x+y \geq k \in 2\Z$,  then $e^\pm(x) +
e^\mp(y) \geq k$, thus (3) holds since $d_{i-1}/m + d_i/m \geq 2r_i$.
Suppose that $d'_{i-1} = e^+(d_{i-1}/m)$ and so $d'_i =
e^-(d_i/m)$. We show that (1) holds.  We have that
$$2r_i \geq d_{i-1}/m - d_i/m
\geq d'_{i-1} - d'_i - 2,$$ so if (1) fails then $d'_{i-1} - d'_{i}
= 2r_i + 2$.   But $d'_{i-1} + d'_i \geq 2r_i$ since we have shown (3)
already, and hence we get that $d'_{i-1} \geq 2r_i + 1$, a contradiction
with $d'_{i-1} \leq 2r_i$. The other cases are similar.
\end{proof}

\begin{theorem}\label{toric_generators}
The toric ring $\C[S_\br]$ is generated
by elements of degrees one and two; furthermore, $(R_{2\br})_0$ is
generated by elements of degree one.
\end{theorem}

\begin{proof}
The first statement is a direct consequence of Lemmas
\ref{odd_degree} and \ref{even_degree}.  The second statement
follows from Lemma \ref{even_degree}.
\end{proof}

\subsection{The word problem for $S_\br$ and the relations for $\C[S_\br]$.}
\hfill

\medskip

We will actually solve the presentation problem for $\C[S_\br]$ by solving
the seemingly more difficult \emph{word problem} for the graded semigroup
$S_\br = \cup_{N \geq 0} D(N)(\Z)$.
Our technique is to
define a normal form for words in $S_\br$ expressed in terms of
degree one and degree two elements, then show that any word
can be brought into normal form by a sequence of quadric relations.

\begin{definition}
Let $\xi_{2m+1} : D(2m+1)(\Z) \to D(1)(\Z)$ be given by
$\xi_{2m+1}(\bd) = \bd'$ where $\bd'$ is as in the proof of Lemma
\ref{odd_degree}.
\end{definition}

\begin{definition}
Let $A$ be an integer matrix.  Let the $j$th column of $A$ be
denoted $c_j(A)$.  If each column of $A$ belongs to either
$D(1)(\Z)$ or $D(2)(\Z)$ then we say that $A$ is a $D$-matrix.
($D$-matrices represent monomials in $\C[S_\br]$ which are products of
degree one and degree two generators.) The elements of $D(1)(\Z)$ (resp.\ $D(2)(\Z)$) are said to
have degree one (resp.\ two). We define $\deg(A)$ to be the sum of
the degrees of the columns of $A$ whenever $A$ is a $D$-matrix.
\end{definition}

\begin{definition}(normal form)
Suppose that $A$ is a $D$-matrix. Let
$$\ba = (a_1,\ldots,a_{n-1})
= \sum_j c_j(A) \in D(\deg(A))(\Z).$$

Suppose that $\deg(A) = 2m$ is even.   Let $\mathcal{J}_\ba =
\{i_1,\ldots,i_k\}$ be the set of all $i$, $2 \leq i \leq n-1$,
such that $a_{i-1} \leq 2m r_i$ and $a_i \leq 2m r_i$, where $i_t <
i_{t+1}$ for all $t$, $1 \leq t < k$. Let $i_0 = 1$ and let
$i_{k+1} = n$. (Note that $\mathcal{J}_\ba = \mathcal{J}_\ba^1$ as
in the proof of Lemma~\ref{even_degree}.) Suppose $A$ has the
following properties:
\begin{itemize}
\item[(N0)] Each column of $A$ has degree two. \item[(N1)] For
each $i$ the row entries $a_{i,j}$ satisfy $|a_{i,j} - a_i/m| <
2$. \item[(N2)] For $i_{2t} \leq i < i_{2t+1}$, row $i$ is weakly
increasing. For $i_{2t+1} \leq i < i_{2t+2}$, row $i$
is weakly decreasing.
\end{itemize}
Then we say that $A$ is in normal form.

Now suppose that $\deg(A) = 2m+1$ is odd.  Then we say that $A$ is
in normal form if the first column of $A$ is equal to
$\xi_{2m+1}(\ba)$ and if the matrix $A'$ obtained from $A$ by
removing the first column is in normal form.
\end{definition}

\begin{example}\label{normalformexample1}
Here is an example of two $D$-matrices $A$ and $B$ where $B$ is the normal form 
representative of $A$.  This for the case $\br = (1,1,1,1,1,1,1,1)$.   
$$A = 
\begin{pmatrix}
1 & 2 & 2 \\
2 & 2 & 0 \\
1 & 0 & 2 \\
2 & 2 & 4 \\
3 & 0 & 4 \\
2 & 2 & 2 \\
1 & 2 & 2 \\ 
\end{pmatrix} \qquad B = 
\begin{pmatrix}
1 & 2 & 2 \\
0 & 2 & 2 \\
1 & 0 & 2 \\
2 & 2 & 4 \\
1 & 2 & 4 \\
2 & 2 & 2 \\
1 & 2 & 2 \\
\end{pmatrix}
$$
The lattice point given by $A$ and $B$ is $\bd = (5,4,3,8,7,6,5) \in D(5)(\Z)$; this 
is the sum of the columns.  According to the definition of normal form, the 
first column should be equal to  $\xi_5(\bd) = (1,0,1,2,1,2,1) \in D(1)(\Z)$; this is indeed 
the first column of $B$.  Now $\bd' = \bd - \xi_5(\bd) = (4,4,2,6,6,4,4) \in D(4)(\Z)$.
Note that $\mathcal{J}_{\bd'} = \{2,3,7\}$.  So the two remaining degree two columns 
of $B$ should have the property that rows $1,3,4,5,6$ weakly increase and rows 
$2,7$ weakly decrease.  Furthermore, the entries in a given row shouldn't differ by more 
than $2$.  This is true for the last two columns of $B$, and so $B$ is in normal form.  
Hence $B$ the normal representative among the set of all tuples of degree one and degree 
two lattice points whose components sum to $\bd = (5,4,3,8,7,6,5)$.
\end{example}

\begin{lemma}(uniqueness)
For any $\ba = (a_1,\ldots,a_{n-1})\in D(\ell)(\Z)$ there is at
most one matrix $A$ in normal form such that the columns of $A$
sum to $\ba$.
\end{lemma}

\begin{proof}
Suppose $\ell = 2m = \deg(A)$ is even and $A$ is in normal form.
Then each column of $A$ is degree two so all the matrix entries
are even integers and there are $m$ columns. For each $i$ let
$k_i$ be an even integer such that $k_i \leq a_i/m \leq k_i + 2$.
By condition (N1) we know that each $a_{i,j}$ is either $k_i$ or
$k_i + 2$. Let $t_i$ be the number of $a_{i,j}$ equal to $k_i$.
Then, $t_i k_i + (m-t_i)(k_i+2) = a_i$, so $2t_i = m (k_i+2) -
a_i$, and thus $t_i$ is determined by the value of $a_i$. Finally
the monotonicity condition (N2) determines each
$a_{i,j}$.

Suppose $\ell = 2m+1 = \deg(A)$ is odd and $A$ is in normal form.
The first column of $A$ must be equal to $\xi_{2m+1}(\ba)$ so it
is determined.  Now the matrix $A'$ which is $A$ with the first
column removed is degree $2m$ and is in normal form, so its 
entries are determined by the argument given above for matrices of
even degree.
\end{proof}

\begin{definition}
We say that two $D$-matrices $A$ and $B$ are equivalent if the sum
of the columns of $A$ is equal to the sum of the columns of $B$.
(Note that each equivalence class contains at most one
representative in normal form by the above lemma --- 
later we will
see that a normal form representative always exists.)
\end{definition}

\begin{definition}
Let $A$ be a $D$-matrix and let $\bd_1$ and $\bd_2$ be two
different columns of $A$.  We define operations of types (F2),
(F3), and (F4) as follows.
\begin{itemize}
\item[(F2)] If $\deg(\bd_1) = \deg(\bd_2) = 1$ then remove columns
$\bd_1$ and $\bd_2$ and place $\bd_1 + \bd_2$ as the last column.
\item[(F3)] If $\deg(\bd_1) = 1$ and $\deg(\bd_2) = 2$, let $\bd =
\bd_1 + \bd_2$.  If $\bd_1$ precedes $\bd_2$ then replace $\bd_1$ with
$\xi_3(\bd)$ and replace $\bd_2$ with $\bd - \xi_3(\bd)$.
If $\bd_2$ precedes $\bd_1$ then replace $\bd_2$ with
$\xi_3(\bd)$ and replace $\bd_1$ with $\bd - \xi_3(\bd)$.
\item[(F4)] If $\deg(\bd_1) =
\deg(\bd_2) = 2$ then choose $\bd'_1$ and $\bd'_2$ each of degree
two such that $\bd'_1 + \bd'_2 = \bd_1 + \bd_2$.  Replace
$\bd_1$ with $\bd'_1$ and replace $\bd_2$ with $\bd'_2$.
\end{itemize}

\end{definition}

\begin{lemma}\label{relations_lemma}
Suppose that $A$ is a $D$-matrix. Then there is a finite sequence
$(A_0,A_1,\ldots,A_p)$ of equivalent $D$-matrices where $A_0 = A$,
the final matrix $A_p$ is in normal form, and for each $i$ the matrix
$A_{i+1}$ is obtained from $A_i$ by a single operation of type (F2), (F3), or
(F4).  (In the special case that $A$ is even degree and all columns are
degree two then all these operations
are of type (F4).)
\end{lemma}

\begin{proof}
First note that (F2) operations can be applied to any pair of degree one columns
until either every column is degree two (when $\deg(A)$ is even)
or there is only one column of degree one (when $\deg(A)$ is odd.)
Assume now that $A$ has at most one column of degree one.

Suppose that $\deg(A) = 2m$ is even, so that each column of $A$ is degree two
and there are $m$ columns.
For $\bd = (d_1,\ldots,d_{n-1}) \in D(4)(\Z)$, let
$\mathcal{J}_\bd = \mathcal{J}_\bd^0$ where $\mathcal{J}_\bd^0$ is
as in the proof of Lemma~\ref{even_degree},
and let $f^-,f^+ : D(4)(\Z) \to D(2)(\Z)$
be given by $f^-(\bd) = \bd'$ and $f^+(\bd) = \bd''$ where
$\bd'$ and $\bd''$ are again as in the proof of Lemma~\ref{even_degree}.
Let an $f^-,f^+$ operation be the following.  Choose $j,j'$ such
that $1 \leq j < j' \leq m$.  Replace columns $c_j(A), c_{j'}(A)$
with $f^-(c_j(A) + c_{j'}(A)),f^+(c_j(A) + c_{j'}(A))$. Then,
$$c_j(A) + c_{j'}(A) = f^-(c_j(A)+c_{j'}(A)) + f^+(c_j(A) +
c_{j'}(A)),$$ so an $f^-,f^+$ operation is of type (F4).
Let $\ba = (a_1,\ldots,a_{n-1}) = \sum_j c_j(A)$.
We claim that
after a finite number of $f^-,f^+$ operations, each entry
$a_{i,j}$ of row $i$ satisfies $|a_{i,j} - a_i/m| < 2$.
The $i$th row (which has even integer entries that sum to $a_i$)
is of minimal distance (using the standard Euclidean metric)
from the constant
vector $(a_i/m,\ldots,a_i/m)$ iff
$|a_{i,j}-a_i/m| < 2$ for each $j$.
Suppose $x$ and $y$ are even integers.  It is easy to check that
$$(x-a_i/m)^2+(y-a_i/m)^2 \geq
\Big(e^+\Big(\frac{x+y}{2}\Big)-a_i/m\Big)^2 +
\Big(e^-\Big(\frac{x+y}{2}\Big)-a_i/m\Big)^2$$ and the inequality
is strict iff $|x-y| \geq 4$. Hence $f^-,f^+$ operations cannot
take the $i$th row further from the constant vector
$(a_i/m,\ldots,a_i/m)$. Now suppose that the $i$th row is as
close as possible to $(a_i/m,\ldots,a_i/m)$ by applying $f^-,f^+$
operations. Suppose there is some $a_{i,j}$ such that $|a_{i,j} -
a_i/m| \geq 2$.  Then there is some $j'$ such that $|a_{i,j'} -
a_{i,j}| > 2$ since $\sum_j a_{i,j} = a_i$.  Since $a_{i,j'} -
a_{i,j}$ is even we have $|a_{i,j'} - a_{i,j}| \geq 4$. But now an
$f^-,f^+$ operation on columns $j,j'$ places row $i$ strictly
closer to $(a_i/m,\ldots,a_i/m)$, a contradiction. Therefore after
sufficiently many $f^-,f^+$ operations the resulting matrix
satisfies (N1).  Assume now that $A$ satisfies (N1). We shall now
switch to a different kind of (F4) operation (which does not
disrupt (N1)) which will eventually give us a matrix that also
satisfies (N2). Recall the definition of $\mathcal{J}_\ba$.  We
have $\mathcal{J}_\ba = \{i_1,\ldots,i_k\}$ is the set of all $i$,
$2 \leq i \leq n-1$, such that $a_{i-1} \leq 2m r_i$ and $a_i \leq
2m r_i$, where $i_t < i_{t+1}$ for all $t$, $1 \leq t < k$. Let
$i_0 = 1$ and let $i_{k+1} = n$. Let $$D(4)(\Z)_\ba = \{\bd \in
D(4)(\Z) \mid \mathcal{J}_\bd^0 \subset \mathcal{J}_\ba \subset
\mathcal{J}_\bd^1\}.$$ Let 
$$g^- : D(4)(\Z)_\ba \to D(2)(\Z),$$
$$g^-(\bd)_i =
\left \{ \begin{array}{ll}
e^-(d_i/2) & \textrm{for $i_{2t} \leq i < i_{2t+1}$}, 2t \leq k,\\
e^+(d_i/2) & \textrm{for $i_{2t+1} \leq i < i_{2t+2}$}, 2t+1
\leq k.
\end{array} \right .$$
Let
$$g^+ : D(4)(\Z)_\ba \to D(2)(\Z),$$
$$g^+(\bd)_i =
\left \{ \begin{array}{ll}
e^+(d_i/2) & \textrm{for $i_{2t} \leq i < i_{2t+1}$}, 2t \leq k,\\
e^-(d_i/2) & \textrm{for $i_{2t+1} \leq i < i_{2t+2}$}, 2t+1.
\leq k
\end{array} \right .$$

We claim that for any two columns $c_j(A)$, $c_{j'}(A)$
the sum $\bd = (d_1,\ldots,d_{n-1}) = c_j(A) + c_{j'}(A)$
is a member of $D(4)(\Z)_\ba$.
First we show that $\mathcal{J}_\ba \subset \mathcal{J}_\bd^1$.
Suppose that $i \in\mathcal{J}_\ba$.
Then $a_{i-1}/m \leq 2r_i$ and $a_i/m \leq 2r_i$, so the entries in
rows $i-1$ and $i$ are at most $2r_i$ since $|a_{i-1,j} - a_{i-1}/m| < 2$
and $|a_{i,j} - a_i/m| < 2$ for all $j$.
Hence the sum of any two entries in row $i-1$
is at most $4r_i$ and the sum of any two entries in row $i$ is at
most $4r_i$.  Therefore each of $d_{i-1}$ and $d_i$ is at most $4r_i$ so
$i \in \mathcal{J}_\bd^1$.
Next we show that $\mathcal{J}_\bd^0 \subset \mathcal{J}_\ba$.
Suppose $i \in \mathcal{J}_\bd^0$.  This means that $d_{i-1} + d_i = 4r_i$ and
each of $d_{i-1}/2$ and $d_i/2$ is an odd integer.  Since $a_{i-1,j} + a_{i,j} \geq 2r_i$
and $a_{i-1,j'} + a_{i,j'} \geq 2r_i$ and
$$a_{i-1,j} + a_{i-1,j'} + a_{i,j} + a_{i,j'} = d_{i-1} + d_i = 4r_i,$$
we have that $a_{i-1,j} + a_{i,j} = 2r_i$ and $a_{i-1,j'} + a_{i,j'} = 2r_i$.
Suppose by way of contradiction that $a_{i-1} > 2mr_i$.  Then each entry of row
$i-1$ is at least $2r_i$.  Then $a_{i-1,j} = a_{i-1,j'} = 2r_i$ and
$a_{i,j} = a_{i,j'} = 0$.  But now $d_i = a_{i,j} + a_{i,j'} = 0$ which contradicts
that $d_i/2$ is an odd integer. Therefore $a_{i-1} \leq 2mr_i$.  Similarly we can show
$a_i \leq 2mr_i$.  Hence $i \in \mathcal{J}_\ba$.

We define a $g^-,g^+$ operation to be the following.
Let $j < j'$ and replace columns $c_j(A)$, $c_{j'}(A)$ with
$g^-(c_j(A)+c_{j'}(A))$, $g^+(c_j(A) + c_{j'}(A))$ in that order.
Clearly any such $g^-,g^+$ operation is of type (F4) and it
preserves the inequalities $|a_{i,j}-a_i/m| < 2$.  We claim that a
finite number of such operations results in a matrix in normal
form. First notice that $g^-,g^+$ operations don't change the
multi-set of entries in any given row since they preserve the (N1) condition.
We determine how
$g^-,g^+$ operations affect the order of the row entries. The
output of $g^-$ and $g^+$ is determined by the type of interval
$i$ belongs to; either $i_{2t} \leq i < i_{2t+1}$ for some $t$ or
$i_{2t+1} \leq i < i_{2t+2}$ for some $t$.  Let us examine the
case $i_{2t} \leq i < i_{2t+1}$.  Here $g^-$ applies the $e^-$
rule and $g^+$ applies the $e^+$ rule.  Hence the result of a
$g^-,g^+$ operation to columns $j,j'$ with $j<j'$ is to put
entries $a_{i,j},a_{i,j'}$ into (weakly) increasing order.  After
applying these operations to all pairs $j,j'$, the resulting
$i$th row is weakly increasing. The case $i_{2t+1} \leq i <
i_{2t+2}$ is similar; this row will be weakly decreasing after
$g^-,g^+$ operations are performed on all pairs of columns.

Suppose $\deg(A) = 2m+1$ is odd, so there is one column of degree one and
$m$ columns of degree two.
Apply a single (F3) operation so that the first column is
the degree one column and columns $2$ through $m+1$ are degree two.
Always let $A'$ denote $A$ without the first column.
We will show after enough operations of types (F3) and (F4)
that the first column is $\xi_{2m+1}(\ba)$ and that $A'$ satisfies conditions
(N0) and (N1) for normality.  Then $g^-,g^+$ operations can be performed on
$A'$ so that $A'$ will eventually satisfy (N2).

The $i$th row must satisfy that $a_{i,1} \equiv (r_1 + \cdots +
r_i) \mod 2$, each $a_{i,j}$ is even for $j \geq 2$, and the sum
$\sum_j a_{i,j} = a_i$. Clearly row $i$ is closest to the vector
$$v_i = \Bigg(\frac{a_i}{2m+1},\frac{2a_i}{2m+1},\frac{2a_i}{2m+1},\ldots,\frac{2a_i}{2m+1}\Bigg) \in \Z^{m+1}$$
iff
$$(*) \quad |a_{i,1} - a_i/(2m+1)| < 1 \; \mbox{ and } \; |a_{i,j} - 2a_i/(2m+1)| < 2
\mbox{ for all } j \geq 2.$$  These inequalities are necessary
for the first column $c_1(A)$ to
be $\xi_{2m+1}(\ba)$ and for $A'$ to satisfy (N1).  If each row satisfies ($\ast$) then
in fact $c_1(A) = \xi_{2m+1}(\ba)$ and $A'$ satisfies (N1).

Suppose that ($\ast$) holds for row $i$.
We claim that
operations of types (F3) and (F4) preserve ($\ast$).
We have that  $|2a_{i,1} - a_{i,j}| < 4$ for each $j \geq 2$.
But $2a_{i,1} - a_{i,j}$ is even so in fact $|2a_{i,1} - a_{i,j}| \leq 2 < 3$.  Therefore,
$|(a_{i,1} + a_{i,j})/3 - a_{i,1}| < 1$.  But this implies that $\xi_3(c_1(A) + c_j(A))_i =
a_{i,1}$ since $a_{i,1}$ has parity $(r_1 + \cdots + r_i) \mod 2$ and is less than one unit from
$(a_{i,1}+a_{i,j})/3$.  Therefore row $i$ is fixed by any (F3) operation. On the other hand
if an (F4) operation is applied to columns $j$ and $j'$ then it either fixes $a_{i,j}$ and
$a_{i,j'}$ or swaps their order since $|a_{i,j} - a_{i,j'}| \leq 2$.

Suppose that row $i$ is as close as possible to $v_i$ by applying (F3) and (F4) operations.
Suppose by way of contradiction that $|a_{i,1} - a_i/(2m+1)| \geq 1$.  Then there is some
$j_0 \geq 2$ such that $a_i/(2m+1)$ is strictly between $a_{i,j_0}/2$ and $a_{i,1}$ since $a_i/(2m+1)$ is the
weighted average of the entries in row $i$, where $a_{i,1}$ is weighted by $1$ and $a_{i,j}$ is
weighted by $2$ for each $j \geq 2$.  Therefore $|2a_{i,1} - a_{i,j_0}| > 2$.  But
$2a_{i,1} - a_{i,j_0}$ is even so in fact $|2a_{i,1} - a_{i,j_0}| \geq 4$.   So
$|(a_{i,1} + a_{i,j_0})/3 - a_{i,1}| \geq 4/3$.  Without loss of generality suppose that
$a_{i,1} < a_{i,j_0}/2$.  Then we have 
$$a_{i,1} < (a_{i,1} + a_{i,j_0})/3 < a_{i,j_0}/2.$$  Let $k$ be the nearest integer of parity
$(r_1 + \cdots + r_i) \mod 2$ to $(a_{i,1} - a_{i,j_0})/3$.  Then we have $a_{i,1} < k \leq a_{i,j_0}/2$.
Let $\delta_i$ be the change in the distance between row $i$ and $v_i$ after applying an (F3) operation
to columns $1$ and $j_0$.  Let $a = a_i/(2m+1)$ and let $t = k - a_{i,1}$.
Then
$$\delta_i = (a_{i,1} + t - a)^2 + (a_{i,j_0} - t - 2a)^2 - (a_{i,1}-a)^2 - (a_{i,j_0}-2a)^2
= 2t(t - (a_{i,j_0} - a_{i,1} - a)).$$
But we know that $0 < t \leq a_{i,j_0}/2 - a_{i,1} < a_{i,j_0} - a_{i,1} - a$.  The first inequality follows
from the fact that $a_{i,1} < k$ and
the last inequality follows from that fact that $a_{i,j_0} > 2a = 2a_i/(2m+1)$.  Hence $\delta_i$ is negative
which means an (F3) operation takes row $i$ strictly closer to $v_i$, a contradiction.  Hence
$|a_{i,1} - a_i/(2m+1)| < 1$.  Now by our argument above for even degree matrices, we must have that
the remaining entries $a_{i,j}$ differ by at most $2$ from one another ((F4) operations can accomplish this)
and consequently we also have that $|a_{i,j} - 2a_i/(2m+1)| < 2$ for each $j \geq 2$.  Therefore, working row by
row, we end up with a matrix $A$ such that $c_1(A) = \xi_{2m+1}(\ba)$ and $A'$ satisfies (N1).  Now apply
$g^-,g^+$ operations to $A'$ so that finally $A'$ satisfies (N2) as well.
\end{proof}

\begin{corollary}
For any $D$-matrix $A$,
there is a unique matrix $\mathcal{N}(A)$ in normal form which is
equivalent to $A$.
\end{corollary}

\begin{theorem}\label{toric_relations}
The ideal of $\C[S_\br] \cong (R_\br)_0 $ is
generated by quadratic relations of degrees two, three, and four.
Furthermore, the ideal of $(R_{2\br})_0$ is generated by quadratic relations.
\end{theorem}

\begin{proof}
One only needs to determine if two monomials in degree one and two variables
are equal.  This corresponds to deciding if two $D$-matrices are equivalent, which
is true iff they have the same normal form.
Operations of types (F2),(F3), and (F4) correspond to degree two, degree three, and
degree four relations in the ideal of $\C[S_\br]$.   By Lemma~\ref{relations_lemma} these operations
are enough to place any $D$-matrix $A$ into its normal form $\mathcal{N}(A)$.  Hence
relations up to degree four must generate the ideal of $\C[S_\br]$.  For the case
of $(R_{2\br})_0$ a $D$-matrix is of even degree and
we only need type (F4) operations to place it into normal form.  In this case an (F4) operation
corresponds to a quadratic relation.
\end{proof}

\section{Lifting generators and relations from the toric
fiber}\label{lifting}

We wish to get a presentation for our ring $R_\br$ by lifting the presentation 
for the associated graded ring $\gr(R_\br)$ that we found in \S \ref{toricring}.  
It is a basic fact that generators of an 
associated graded \emph{algebra} or a \emph{module} may be lifted to generators 
of the original object (whether it be an algebra or module.)  Here we will use both facts, 
that is, we will lift generators of the algebra $\gr(R_\br)$ to get generators of $R_\br$, and we 
will lift generators of the ideal of $\gr(R_\br)$ (a module) to get generators of the ideal of $R_\br$.

We begin with the following lemma which will be basic in what
follows. We leave the proof to the reader. For background on
filtrations and gradings we refer to \cite{Bourbakicomalg}.
Throughout this section all filtrations will be increasing and
indexed by the nonnegative integers $\Z^+$.

\begin{lemma}\label{generationlemma}
Suppose that $M$ is a filtered module over a filtered ring $R$ and
that their filtrations are compatible in the sense that
$$F_i(R) \otimes F_j(M) \subset F_{i+j}(M).$$
Suppose that $x_1,x_2,\ldots,x_n$ are
elements of $M$ such that their images $\overline{x}_i, 1 \leq i
\leq n$, under the leading term map generate $\gr(M)$ as a $\gr(R)$ module.
Then $x_i, 1 \leq i
\leq n$ generate $M$.
\end{lemma}

\begin{remark} An analogous argument shows that if the images in $\gr(R)$ of
a finite set of elements $r_1, r_2,\ldots, r_n$ of $R$ generate $\gr(R)$ then the
elements $r_1$. $r_2$, \ldots, $r_n$ generate $R$.
\end{remark}

Our goal in this section is to prove the statement for {\it
relations} that is the analogue of the statement in the remark for
generators.

\begin{definition}
Let $M$ be a filtered module and $x \in M$. We define the filtration level
(or order) $v(x) \in \Z^+$ of $x$ to be the smallest $n$ such that $x \in F_n(M)$.
\end{definition}

Assume that $R$ is graded as a vector space and that we have
chosen homogeneous generators $f_1,f_2,\ldots,f_n$ for $R$  such
that the images
$\overline{f}_1,\overline{f}_2,\ldots,\overline{f}_n$ of these
generators in $\gr(R)$ generate $\gr(R)$. We assume the degree of
$f_i$ is $e_i, 1 \leq i \leq n$.

We obtain two exact sequences.
$$
\begin{CD}
I @ >\iota >> \C[x_1,x_2,\cdots,x_n] @ > \pi >> R
\end{CD}
$$
and
$$
\begin{CD}
J @ >>> \gr(\C[x_1,x_2,\cdots,x_n]) @ > \pi >> \gr(R).
\end{CD}
$$
Here $\pi$ sends $x_i$ to $f_i, 1 \leq i \leq n$. In the above
the polynomial ring $\C[x_1,x_2,\cdots,x_n]$ is a weighted polynomial
ring, the variable $x_i$ has weight $e_i$. We define a filtration
on $R$ by defining the filtration level of $r$ to be the minimum
of the degrees of the polynomials in $\pi^{-1}(r)$. The reader will
verify that this filtration coincides with the quotient filtration of the
standard filtration on $\C[x_1,x_2,\cdots,x_n]$. We remind the
reader that the quotient filtration is characterized by the fact that
the induced map on each filtration level is a surjection,
see \cite[pg.~164]{Bourbakicomalg}.

We note that
$\gr(\C[x_1,x_2,\cdots,x_n]))$ is the polynomial ring
$\C[\bar{x}_1,\bar{x}_2,\cdots,\bar{x}_n]$.
We leave the reader the task of proving (by induction on
the filtration level):

\begin{lemma}
Suppose $R$ is a filtered $\C$--algebra which is graded as a vector space
and $f_1,\ldots,f_n$ have the property that their images $\overline{f}_1,
\ldots, \overline{f}_n$ generate $\gr(R)$. Then the given filtration on
$R$ coincides with the quotient
filtration associated to the surjection $\pi:\C[x_1,x_2,\cdots,x_n] \to R$
given by $\pi(x_i) = f_i$.
\end{lemma}

\begin{example}
We give the following example to show what can go wrong if
$\overline{f}_1,\ldots,\overline{f}_n$ do not generate $\gr(R)$. We give
this example because a similar phenomenon occurs for the case of
equilateral hexagons. Consider the affine coordinate ring $R$ of the
saddle surface $z = xy$. Then $R$ is generated by $x$ and $y$.
We give $R$ the filtration that is the quotient of the filtration
on $\C[x,y,z]$, so the images $\overline{x}, \overline{y}$ and
$\overline{z}$ have degree  one in $R$. We also
have a surjection $\pi:\C[u,v] \to R$ given by $\pi(u) = x$ and
$\pi(v) = y$. Clearly $\pi$ is not onto at filtration level one.
There is no contradiction with the lemma above because  the images
$\overline{x}$ and $\overline{y}$ do not generate $\gr(R)$.
\end{example}

Since we give $I$ the filtration induced as a submodule of the
polynomial ring both $I$ and $R$ have the filtrations needed to
apply \cite[pg.~169, Prop.~2]{Bourbakicomalg}  to deduce
that we have an exact sequence
$$
\begin{CD}
\gr(I) @ >\gr(\iota)>> \gr(\C[x_1,x_2,\cdots,x_n]) @ >\gr(\pi)>>
\gr(R).
\end{CD}
$$
and consequently $\gr(\iota): \gr(I) \to J$ is an isomorphism.

We are now ready to state and prove the result we want on lifting
relations from $\gr(R)$ to $R$. We emphasize that we are assuming
that the generators for $R$ map to generators for $\gr(R)$ under
the leading term map.

\begin{proposition}\label{liftingrelations}
Suppose $p_1,p_2,\ldots,p_k \in \gr(\C[x_1,x_2,\dots,x_n])$
generate the ideal of relations in the given generators for
$\gr(R)$. Then
\begin{enumerate}
\item There exist lifts $\tilde{p}_i, 1 \leq i \leq k$, to
$C[x_1,x_2,\ldots,x_n]$ such that for all $i$ the polynomial
$\tilde{p}_i$ is a relation for $R$. \item For any choice of such
lifts $\tilde{p}_i, 1 \leq i \leq k$ the lifts generate the ideal
of relations of $R$.
\end{enumerate}
\end{proposition}

\begin{proof}
Since we have shown that $J \equiv \gr(I)$ the first statement in
the proposition is obvious (since the leading term map is onto by
definition of $\gr(I)$). However the lift of a homogeneous element
will usually not be homogeneous (the ideal $I$ may not  contain
any nonzero homogeneous elements). The second statement follows
from Lemma~\ref{generationlemma} --- the images of the lifts
generate the ideal $\gr(I)$ so the lifts generate $I$.
\end{proof}

\section{The projective coordinate ring of $M_\br$}\label{ring}

Now we apply the lifting results of the previous section to get a
presentation of the ring $R_\br$.  We begin by carrying over some
definitions from the toric ring $(R_\br)_0$ to $R_\br$.  Recall
the definition of normal form for $D$-matrices.

\begin{definition}
For any $\bd \in D(m \br)(\Z)$ let $\tau_\bd \in R_\br$ be the
unique tableau with multi-weight $\Psi(\Phi^{-1}(\bd)) \in SS(ma
\varpi_2, m\br)$. If a $D$-matrix is in normal form with columns
$\bd_1,\ldots,\bd_s$ we shall say the product $\tau_{\bd_1}
\tau_{\bd_2} \cdots \tau_{\bd_s} \in R_\br$ is a normal monomial.
Also we shall say that $\bar{\tau}_{\bd_1}
\bar{\tau}_{\bd_2} \cdots \bar{\tau}_{\bd_s} \in (R_\br)_0$ is a
normal monomial.
\end{definition}

\begin{example}
We return to Example \ref{normalformexample1}, and see what are the related 
monomials in tableaux.  We had:
$$A = 
\begin{pmatrix}
1 & 2 & 2 \\
2 & 2 & 0 \\
1 & 0 & 2 \\
2 & 2 & 4 \\
3 & 0 & 4 \\
2 & 2 & 2 \\
1 & 2 & 2 \\ 
\end{pmatrix} \qquad B = 
\begin{pmatrix}
1 & 2 & 2 \\
0 & 2 & 2 \\
1 & 0 & 2 \\
2 & 2 & 4 \\
1 & 2 & 4 \\
2 & 2 & 2 \\
1 & 2 & 2 \\
\end{pmatrix}
$$
The associated monomials in tableaux are:
$$
m_A = \tableau{1}{3}\tableau{2}{6}\tableau{4}{7}
\tableau{5}{8} \cdot \tableau{1}{2}\tableau{1}{3}
\tableau{2}{3}\tableau{4}{5}\tableau{4}{5}\tableau{6}{7}
\tableau{6}{8}\tableau{7}{8} \cdot \tableau{1}{2} 
\tableau{1}{2}\tableau{3}{5}\tableau{3}{6}\tableau{4}{6}
\tableau{4}{7}\tableau{5}{8}\tableau{7}{8} \; ,
$$
$$
m_B = \tableau{1}{2}\tableau{3}{5}\tableau{4}{7}
\tableau{6}{8} \cdot \tableau{1}{2}\tableau{1}{3}
\tableau{2}{3}\tableau{4}{5}\tableau{4}{6}\tableau{5}{7}
\tableau{6}{8}\tableau{7}{8} \cdot \tableau{1}{2} 
\tableau{1}{3}\tableau{2}{5}\tableau{3}{6}\tableau{4}{6}
\tableau{4}{7}\tableau{5}{8}\tableau{7}{8} \; .
$$
Here $m_B$ is a normal monomial since $B$ is in normal form, whereas $m_A$ is not
a normal monomial because $A$ is not in normal form.  However, 
$m_B$ and $m_A$ have the same LG-filtration level.
\end{example}

\begin{definition}\label{normalformbasisdefinition}
Suppose $k = 2m$ is even.  Choose tableaux $\tau_{i,j}$ such that
$1 \leq i \leq M_k$, $1 \leq j \leq m$, where $M_k$ is the number
of lattice points in $D(k \br)(\Z)$, and such that
$$\mathcal{N}_{k \br} = \Bigg(\prod_{j=1}^m \tau_{1,j}, \prod_{j=1}^m \tau_{2,j}, \ldots,
\prod_{j=1}^m \tau_{M_k,j}\Bigg) \in (R_\br^{(k)})^{M_k}$$
is the set of normal monomials
of degree $k$,
$$\mathcal{SS}_{k \br} = (\sigma_1,\ldots,\sigma_{M_k})
\in (R_\br^{(k)})^{M_k},
\text{ where}$$
$$\sigma_i = \tau_{i,1} \ast \tau_{i,2} \ast \cdots \ast
\tau_{i,m},$$ such that $\text{LG-deg}(\sigma_i) \leq
\text{LG-deg}(\sigma_{i+1})$ for each $i$, $1 \leq i < M_k$.

Suppose $k = 2m+1$ is odd.  Choose tableaux $\tau_{i,j}$ such that
$1 \leq i \leq M_k$, $1 \leq j \leq m+1$, where $M_k$ is the
number of lattice points in $D(k \br)(\Z)$, and such that
$$\mathcal{N}_{k \br} = \Bigg(\prod_{j=1}^{m+1} \tau_{1,j}, \prod_{j=1}^{m+1} \tau_{2,j}, \ldots,
\prod_{j=1}^{m+1} \tau_{M_k,j}\Bigg) \in (R_\br^{(k)})^{M_k}$$
is the set of normal monomials of degree $k$,
$$\mathcal{SS}_{k \br} = (\sigma_1,\ldots,\sigma_{M_k})
 \in (R_\br^{(k)})^{M_k}, \text{ where}$$
$$\sigma_i = \tau_{i,1} \ast \tau_{i,2} \ast \cdots \ast
\tau_{i,m+1},$$ such that $\text{LG-deg}(\sigma_i) \leq
\text{LG-deg}(\sigma_{i+1})$ for each $i$, $1 \leq i < M_k$.

Note that $M_k$ is the number of semistandard tableaux of weight $k
\br$, which is the dimension of $(R_\br)^{(k)}$.
\end{definition}

\begin{remark}
The choice of the ordering of normal monomials is not unique,
since it is possible for two normal monomials to have the same
total LG-degree, which is $\sum_j \text{LG-deg}(\tau_{i,j})$.
\end{remark}

\begin{proposition}
The components of $\mathcal{SS}_{k \br}$ are exactly the
semistandard tableaux of weight $k \br$.
\end{proposition}

\begin{proof}
Let $\mathcal{N}_{k \br}$ be as above.  The normal monomials
of degree $k$ are in bijection with the integral points of the
polytope $SS(ka\varpi_2, k \br)$, by
$$\prod_j \tau_{i,j} \mapsto \sum_j \wt(\tau_{i,j}) =
\wt(\sigma_i).$$ But the semistandard tableaux of weight $k \br$
are in bijection with the integral points of $SS(ka\varpi_2, k
\br)$ by $\sigma \mapsto \wt(\sigma)$.
\end{proof}

\begin{theorem}\label{normalformbasistheorem}
The tuple $\mathcal{N}_{k \br}$ is a basis for $(R_\br)^{(k)}$.
\end{theorem}

\begin{proof}
Let $c_{i,j} \in \Z$ for $1 \leq i,j \leq M_k$ be given by
$$\prod_\ell \tau_{j,\ell} = \sum_{i = 1}^{M_k} c_{i,j} \;
\sigma_i.$$ 
We claim the matrix $[c_{i,j}]$ is unipotent (or upper-triangular). Since the
product of any two semistandard tableaux $\tau \sigma$ is $\tau
\ast \sigma + \sum_i c_i \rho_i$ where each $c_i \in \Z$, each
$\rho_i$ is semistandard, and $\text{LG-deg}(\rho_i) <
\text{LG-degree}(\tau \ast \sigma)$, one can argue by induction on
$s \geq 2$ that
$$\prod_{\ell = 1}^s
\tau_\ell = (\tau_1 \ast \tau_2 \ast \cdots \ast \tau_s) + \sum_i
c_i \rho_i$$ where each $c_i \in \Z$, each $\rho_i$ is
semistandard with $$\text{LG-deg}(\rho_i) < \text{LG-deg}(\tau_1
\ast \tau_2 \ast \cdots \ast \tau_s).$$  Since the tuple
$\mathcal{SS}_{k \br} = (\sigma_1, \ldots, \sigma_{M_k})$
satisfies $\text{LG-deg}(\sigma_i) \leq
\text{LG-deg}(\sigma_{i+1})$, the matrix $[c_{i,j}]$ is
upper-triangular, and since the first term of the sum is the
concatenated product, we get $1$'s on the diagonal of $[c_{i,j}]$.
Since $\mathcal{SS}_{k \br}$
is a basis of $R_\br^{(k)}$, the theorem follows.
\end{proof}

\begin{definition}
Let $C_{k \br}$ be the (unipotent) basis exchange matrix from
$\mathcal{N}_{k \br}$ to $\mathcal{SS}_{k \br}$.
\end{definition}

\begin{proposition}\label{leadingtermprop}
For any two tableaux $\rho_1,\rho_2$ with
$\rho_1 \rho_2 \in R_\br^{(k)}$, there exists $i_0 \leq M_k$
such that
$$\rho_1 \rho_2 = \prod_{j=1}^s \tau_{i_0,j} + \sum_{i < i_0} c_i
\Bigg(\prod_{j=1}^s \tau_{i,j}\Bigg) \t{, where}$$
$$\mathcal{N}_{k \br} = \Bigg(\prod_{j=1}^s \tau_{1,j}, \prod_{j=1}^s \tau_{2,j}, \ldots,
\prod_{j=1}^s \tau_{M_k,j}\Bigg)$$ and
$\bar{\rho}_1 \bar{\rho}_2 = \prod_{j=1}^s \bar{\tau}_{i_0,j}$
is the leading term of the right hand side with respect to LG-degree.
\end{proposition}

\begin{proof}
We have already established the analogous statement in terms of semistandard
tableaux.  There exists $i_0 \leq M_k$
such that
$$\rho_1 \rho_2 = \sigma_{i_0} + \sum_{i < i_0} c'_i
\sigma_i \t{, where}$$
$$\mathcal{SS}_{k \br} = (\sigma_1,\ldots,\sigma_{M_k})$$ and
$\bar{\rho}_1  \bar{\rho}_2 = \bar{\sigma}_{i_0}
= \prod_{j=1}^s \bar{\tau}_{i_0,j}$ is the leading term of
the right hand side with respect
to LG-degree.
\end{proof}

\begin{theorem}\label{relations} \hfill
\begin{enumerate}
\item The ring $R_\br$ is generated by semistandard tableaux in
degrees one and two.  The relations in these tableaux are
generated by quadratic relations of degrees two, three,
and four.
\item The
ring $R_{2 \br}$ is generated by semistandard tableaux of degree
one.  The relations in these tableaux are generated by quadratic
relations.
\end{enumerate}
\end{theorem}

\begin{proof}
The tableaux $\bar{\tau}$ of degrees one and two generate $(R_\br)_0$ so
their lifts $\tau$ must generate $R_\br$.

Now we choose lifts of the relations in $(R_\br)_0$ of types (F2),
(F3), and (F4) to get degree two, three, four relations in the
tableaux of degree one and two.

First consider a relation of type (F2) which is $\bar{\sigma}_1\bar{\sigma}_2
= \bar{\tau}$
where $\sigma_1 \ast \sigma_2 = \tau$ and $\sigma_1$, $\sigma_2$ are degree one.
Write  
$$\sigma_1 \sigma_2 = \sum_{i=1}^{M_2} c_i
\tau_i \quad \quad \quad (F2')$$ where $M_2 = \dim R_\br^{(2)}$,
$(\tau_1,\ldots,\tau_{M_2}) = \mathcal{SS}_{2\br} = \mathcal{N}_{2 \br}$.
There is
some $k$ such that $\tau = \tau_k$, and furthermore
$\tau = \tau_k$ is the leading term of the right hand side
with respect to LG-degree.  Hence this relation is a lift of the
(F2) relation $\bar{\sigma}_1\bar{\sigma}_2 = \bar{\tau}$.

Now consider a relation of type (F3), which is $\bar{\sigma}\bar{\tau} =
\bar{\sigma}' \bar{\tau}'$ where the right hand side is a normal monomial, and
$\bar{\sigma}, \bar{\sigma}'$ are degree one and $\bar{\tau}, \bar{\tau}'$
are degree two.
Write
$$\sigma \tau =
\sum_{i=1}^{M_3} c_i \sigma_i \tau_i \quad \quad \quad (F3')$$
where $M_3 = \dim R_\br^{(3)}$ and $\mathcal{N}_{3 \br} = (\sigma_1 \tau_1,
\ldots, \sigma_{M_3} \tau_{M_3})$.  The normal monomial
$\bar{\sigma}' \bar{\tau}'$ is the
leading term of the right hand side with respect to LG-degree, so
this relation is a lift of the (F3) relation $\bar{\sigma}\bar{\tau} =
\bar{\sigma}' \bar{\tau}'$.

Finally consider a relation of type (F4), which is $\bar{\tau}_1\bar{\tau}_2 =
\bar{\tau}'_1 \bar{\tau}'_2$ where the right hand side is a normal monomial, and
each of $\bar{\tau}_1, \bar{\tau}_2, \bar{\tau}'_1, \bar{\tau}'_2$
is degree two.
Write
$$\tau_1 \tau_2 =
\sum_{i=1}^{M_3} c_i \tau_{i,1} \tau_{i,2} \quad \quad \quad  (F4')$$
where $M_4 = \dim R_\br^{(4)}$ and $\mathcal{N}_{4 \br} = (\tau_{1,1} \tau_{1,2},
\ldots, \tau_{M_4,1} \tau_{M_4,2})$.  The normal monomial
$\bar{\tau}'_1 \bar{\tau}'_2$ is the
leading term of the right hand side with respect to LG-degree, so
this relation is a lift of the (F4) relation $\bar{\tau}_1\bar{\tau}_2 =
\bar{\tau}'_1 \bar{\tau}'_2$.

It follows from Proposition~\ref{liftingrelations} and Theorem \ref{toric_relations} that the
relations in the degree one and degree two tableaux are generated
by the lifted relations above,
$$\rho_1 \rho_2 = \sum_{i=1}^{M_j} c_i
\mathbf{n}_i \quad \quad \quad  (Fj')$$ where $\mathcal{N}_{j \br} =
(\mathbf{n}_1,\ldots,\mathbf{n}_{M_j})$ and $\rho_1 \rho_2 \in
R_\br^{(j)}$, $j = 2,3,4$.

Similarly the semistandard Young tableaux of weight $2\br$
generate $R_{2 \br}$ and the relations in these tableaux are generated by
quadratic relations, which are lifts of relations of type (F4) for
the toric ring $(R_\br)_0$.
\end{proof}

Now recall Theorem~\ref{generatorstheorem} which states that in fact
$R_\br$ is generated by elements of degree one provided that
$|\br| =  \sum_i r_i$ is even.  We shall now improve this theorem to include
a statement about the relations.

\begin{theorem}\label{maintheorem}
Assume that $|\br| = \sum_{i=1}^n r_i$ is even.  Then $R_\br$ is generated by 
the semistandard
tableaux of weight $\br$ and the relations amongst these tableaux are
generated by relations of degree four and less.
If furthermore each $r_i$ is even then the relations among the tableaux
of weight $\br$ are generated by quadratic relations.
\end{theorem}

\begin{proof}
We have that the degree tableaux of weight $\br$ generate $R_\br$ by Theorem~\ref{generatorstheorem}.

Each tableaux of weight $2 \br$ is a quadratic function of tableaux of weight
$\br$.
For each tableaux $\tau$ of weight $2 \br$, write $\tau$ in terms of
tableaux of weight $\br$:
$$\tau = \sum_i a_i(\tau) \sigma_{i,1} \sigma_{i,2}.$$
Let $f(\tau) = \sum_i a_i(\tau) \sigma_{i,1} \sigma_{i,2}$.
The relations of type (F2'), (F3'), and (F4') in
Theorem~\ref{relations} above may be replaced with relations in tableaux
of weight $\br$ by substituting each degree two tableau $\tau$ with
$f(\tau)$.
\end{proof}

\subsection{An explicit algorithm for listing the relations}

First write each degree two tableau $\tau$ in terms of degree one tableaux
(see Theorem~\ref{generatorstheorem}):
$$\tau = f(\tau) = \sum_i a_i(\tau) \sigma_{i,1} \sigma_{i,2}.$$

For each $j = 2,3,4$, do the following:
\begin{enumerate}
\item List the normal monomials $\mathcal{N}_{j \br} = (\mathbf{n}_1,\ldots,
\mathbf{n}_{M_j})$ and semistandard tableaux
$\mathcal{SS}_{j \br} = (\sigma_1,\ldots,\sigma_{M_j})$ weakly increasing
with respect to LG-degree.
\item Compute the change of basis matrix $C_{j \br}$ using
Pl\"ucker relations.
\item Compute $C^{-1}_{j \br}$ (recall $C_{j \br}$ is
unipotent with integer entries).
\item List all non-normal monomials $\rho_1 \rho_2$ of degree $j$.
\item For each such $\rho_1 \rho_2$ above, use Pl\"ucker relations to compute the coefficients
$c'_{i,\rho_1 \rho_2}$ where
$$\rho_1 \rho_2 = \sum_{i=1}^{M_j} c'_{i,\rho_1\rho_2} \sigma_i.$$

\item Compute $[c_{i,\rho_1\rho_2}]_{i=1}^{M_j} =
C_{j \br}^{-1} \, [c'_{i,\rho_1\rho_2}]_{i=1}^{M_j}$.
Then,
$$\rho_1 \rho_2 = \sum_{i=1}^{M_j} c_{i,\rho_1\rho_2} \mathbf{n}_i.$$

\item Replace each degree two tableau $\tau$ occuring in the relation above
with $f(\tau)$.
\end{enumerate}

\subsection{A minimal set of generating relations for $R_{2 \br}$}

The relations given for $R_\br$ above may not be minimal.  The
prime example is that of $n=6$ with each $r_i = 1$, where a single
cubic relation generates the ideal of relations. However, the
presentation given for $R_{\bs}$ when each $s_i$ is even (say $\bs = 2 \br$) 
turns
out to be a minimal presentation.

\begin{theorem}\label{minimalrelations}
A minimal set of generating relations for $R_{2 \br}$ is given by
$$\rho_1 \rho_2 = \sum_{i=1}^{M_4}
c_{i,\rho_1\rho_2} \sigma_i \tau_i$$ where $(\sigma_1
\tau_1,\ldots,\sigma_{M_4} \tau_{M_4}) = \mathcal{N}_{4 \br}$ and
the products $\rho_1 \rho_2$ range over all non-normal products of
degree two tableaux in $R_\br$ (they are degree one in
$R_{2\br}$).
\end{theorem}

\begin{proof}
Let $I$ be the ideal generated by the above relations, and let $J$
be the ideal generated by all but one of them, say
$p(\rho_1,\rho_2) = \tau_1 \tau_2 - \sum_{i=1}^{M_4} c_{i,\rho_1
\rho_2} \sigma_i \tau_i$. We show that $J$ does not contain
$p(\rho_1,\rho_2)$. Suppose $p(\rho_1,\rho_2) \in J$. Then,
$$p(\rho_1,\rho_2) = \sum_\ell a_\ell p(\alpha_\ell,\beta_\ell)$$
for some coefficients $a_\ell$, where the products $\alpha_\ell
\beta_\ell$ are not equal to $\rho_1 \rho_2$ in the polynomial
ring $\C[\gamma_1,\ldots,\gamma_{M_2}]$ where the $\gamma_i$ range
over the tableaux of weight $2 \br$. But the monomial $\rho_1
\rho_2$ does not even occur on the right hand side of the equation
since it is not one of the $\alpha_\ell \beta_\ell$'s nor is it a
normal monomial. Hence we have a contradiction.
\end{proof}

\subsection{Some examples}

Here we will concentrate on the case of equal weights.
When the $r_i$'s are equal
there is a natural action of the
symmetric group on the moduli space, so that the resulting
quotient is the space of \emph{unordered} points on the line.

\subsubsection{Four points}
Let $\br = (1,1,1,1)$. The moduli space $M_\br$ is simply the
projective line $\CP^1$. The invariant semistandard $2$ by $2$
tableaux are
\begin{center} $X = $
\begin{tabular}{| c | c |}
\hline 1  & 2  \\ \hline 3  & 4  \\ \hline
\end{tabular}
$\quad \quad \quad Y = $
\begin{tabular}{| c | c |}
\hline
    1  &     3  \\ \hline
    2  &     4  \\ \hline
\end{tabular}
\end{center}
These give the embedding into $\CP^1$. The moduli space of the
square is one dimensional over $\C$, and hence there can be no
relations between $X$ and $Y$.  Hence the embedding surjects onto
$\CP^1$.  Identify $\CP^1$ with $\C \cup \{\infty\}$ where $[z,1]
\mapsto z \in \C$, and $[1,0] \mapsto \infty$.  Now let
$z_1,z_2,z_3,z_4 \in \C$.  The condition for semistability is that
no three of these points coincide, and the well known cross ratio
function is given by
$$f(z_1,z_2,z_3,z_4) = \frac{(z_1 - z_3)(z_2 - z_4)}{(z_1
- z_2)(z_3 - z_4)} = (X/Y)([z_1,1],[z_2,1],[z_3,1],[z_4,1]).$$

\subsubsection{Five points}
Let $\br = (2,2,2,2,2)$. The moduli space $M_\br$ (the pentagon
space) is more complicated, it is embedded in $\CP^5$ and
satisfies five quadratic equations. It can be shown that $M_\br$
is $\CP^2$ with four points blown up. The generators are
\begin{center} $A = $
\begin{tabular}{| c | c | c | c | c |}
\hline 1 & 1 & 2 & 2 & 3 \\ \hline 3 & 4 & 4 & 5 & 5 \\ \hline
\end{tabular}
$\qquad B = $
\begin{tabular}{| c | c | c | c | c |}
\hline 1 & 1 & 2 & 2 & 4 \\ \hline 3 & 3 & 4 & 5 & 5 \\ \hline
\end{tabular}
$\qquad C = $
\begin{tabular}{| c | c | c | c | c |}
\hline 1 & 1 & 2 & 3 & 3 \\ \hline 2 & 4 & 4 & 5 & 5 \\ \hline
\end{tabular}
\end{center}
\vskip 12pt
\begin{center}
$D = $
\begin{tabular}{| c | c | c | c | c |}
\hline 1 & 1 & 2 & 3 & 4 \\ \hline 2 & 3 & 4 & 5 & 5 \\ \hline
\end{tabular}
$\qquad E = $
\begin{tabular}{| c | c | c | c | c |}
\hline 1 & 1 & 2 & 4 & 4 \\ \hline 2 & 3 & 3 & 5 & 5 \\ \hline
\end{tabular}
$\qquad F = $
\begin{tabular}{| c | c | c | c | c |}
\hline 1 & 1 & 3 & 3 & 4 \\ \hline 2 & 2 & 4 & 5 & 5 \\ \hline
\end{tabular}
\end{center}

These give an embedding of the moduli space of the pentagon
into $\CP^5$.  There are five non-normal quadratic products,
$BC$, $AE$, $AF$, $BF$, and $CE$.  The generating set of relations
is
$$BC = AD, \qquad
AE = BD - D^2 + DE, \qquad AF = CD - D^2 + DF$$
$$BF = D^2 - DE, \qquad CE = D^2 - DF$$

\subsubsection{Six points}
The space of equilateral hexagons is a cubic hypersurface in $\CP^4$.
We leave this as an exercise for the reader; otherwise
see \cite[pg.~17]{DO} for more details.

\subsubsection{Eight points}
By straightforward computer check, one may verify that the quadratic relations
generate the cubic and quartic relations when $n=8$.  By our main
theorem on relations, we know that \emph{all} relations are
generated by the quadratic relations.  There are fourteen
generators and fourteen relations.  We will use this calculation
in \cite{HowardMillsonSnowdenVakil}.

\vskip 12pt

\centerline{\emph{Generators for the ring of
octagons:}}
$$A =\tableau{1}{5}\tableau{2}{6}\tableau{3}{7}\tableau{4}{8} \quad
B=\tableau{1}{4}\tableau{2}{6}\tableau{3}{7}\tableau{5}{8} \quad
C=\tableau{1}{4}\tableau{2}{5}\tableau{3}{7}\tableau{6}{8}$$
$$D =\tableau{1}{4}\tableau{2}{5}\tableau{3}{6}\tableau{7}{8} \quad
E=\tableau{1}{3}\tableau{2}{6}\tableau{4}{7}\tableau{5}{8} \quad
F=\tableau{1}{3}\tableau{2}{5}\tableau{4}{7}\tableau{6}{8} $$
$$G =\tableau{1}{3}\tableau{2}{5}\tableau{4}{6}\tableau{7}{8} \quad
H=\tableau{1}{3}\tableau{2}{4}\tableau{5}{7}\tableau{6}{8} \quad
I=\tableau{1}{3}\tableau{2}{4}\tableau{5}{6}\tableau{7}{8} $$
$$J =\tableau{1}{2}\tableau{3}{6}\tableau{4}{7}\tableau{5}{8} \quad
K=\tableau{1}{2}\tableau{3}{5}\tableau{4}{7}\tableau{6}{8} \quad
L=\tableau{1}{2}\tableau{3}{5}\tableau{4}{6}\tableau{7}{8} $$
$$M =\tableau{1}{2}\tableau{3}{4}\tableau{5}{7}\tableau{6}{8} \quad
N=\tableau{1}{2}\tableau{3}{4}\tableau{5}{6}\tableau{7}{8} $$
\centerline{\emph{Relations for the ring of octagons:}}
\begin{align}
0 = &-AH  + AI  + AM - AN  + BF - BG  - BK + BL \notag \\
0 = &-BF  + CE \notag \\
0 = &+CL - CN - DK + DM \notag \\
0 = &-FL + GK \notag \\
0 = &+CG  - CI  - DF  + DH + FI  - FN - GH  + GM \notag \\
0 = &-BG + CG  - CI  + DE - DF  + DI  + FI  - FN 
    -GI  + GN \notag \\
0 = &-AH  + BF  - DF  + DH + FI  - FJ + FL + FM \notag \\
    &-2FN- GH  + HJ - HK + HN \notag \\
0 =  &+EL - EN - GJ + IJ \notag \\
0 =  &+EK - EM - FJ + FM - FN + HJ - HK + IK \notag \\
0 = &-AH  + BF  - DF  + DH  + FI  - FJ + FL + FM \notag \\
    &-2FN- GH  + HJ - HK + IM \notag \\
0 = &-AI  + DE  - DF  + DI  + EK - EL - EM + EN \notag \\
    &+FI  - FJ + FL + FM - 2FN- GI  + HJ - HK 
    +IN \notag \\
0 = &-BK + CJ + EK - EM - FJ + FM \notag \\
    &-FN + JM- KM+ KN \notag \\
0 = &-AM + CG  - CI  + CJ - CL + CN - DF  + DH  \notag \\
    &+FI  - FJ + FL + FM - FN - GH  + JM- KM 
    +MN \notag \\
0 = &-AN - BG  - BK + CG - CI  + CJ + DE  - DF \notag \\
    &+DI  + DJ - DK + DN + EK - EM + FI  - FJ \notag \\
    &+FL + FM - 2FN- GI - GJ + JM+ JN- KM 
    +NN \notag
\end{align}

\bigskip 

Benjamin Howard:
Department of Mathematics,
University of Maryland,
College Park, MD 20742, USA,
bhoward@math.umd.edu

\smallskip

John Millson:
Department of Mathematics,
University of Maryland,
College Park, MD 20742, USA,
jjm@math.umd.edu
 
\smallskip

Andrew Snowden:
Department of Mathematics,
Princeton University,
Princeton, NJ 08544, USA,
asnowden@math.princeton.edu

\smallskip

Ravi Vakil:
Department of Mathematics,
Stanford University, 
Stanford, CA 94305-2125, USA,
vakil@math.stanford.edu

\end{document}